\documentclass[preprint,12pt]{elsarticle}

\usepackage{amssymb}
\usepackage{amsmath}

\usepackage{changes}
\graphicspath{{./images/}}

\usepackage[utf8]{inputenc}
\usepackage[T1]{fontenc}
\usepackage{helvet}
\usepackage{etoolbox}
\usepackage{graphicx}
\usepackage{subcaption}
\usepackage{caption}
\usepackage{booktabs}
\usepackage{mathtools}
\usepackage{amsmath}
\usepackage{amssymb}
\usepackage{placeins}
\usepackage{xcolor} 
\usepackage{colortbl}
\usepackage{enumitem}
\usepackage{fancyhdr}
\usepackage{gensymb}
\usepackage{textcomp}
\usepackage{graphicx}
\usepackage{tikz}
\usepackage{multirow}
\usepackage{float}
\usepackage{subcaption}
\usepackage{array}

\usepackage[colorlinks=true, linkcolor=blue, citecolor=blue, urlcolor=blue, pdfborder={0 0 0}]{hyperref}

\def\dd#1#2{\frac{\partial#1}{\partial#2}}

\newcommand{\eps}{\varepsilon}

\newcommand{\dg}{\partial \Omega}

\newcommand{\real}{\mathbb R}

\newcommand{\vecx}{\textbf{x}}
\newcommand{\vecu}{\textbf{u}}

\newcommand{\vecb}{\textbf{b}}
\newcommand{\vecf}{\textbf{f}}

\newcommand{\vecn}{\textbf{n}}

\newcommand{\mb}{\mathcal{B}}

\newcommand{\argmin}[1]{\text{arg}\,\underset{#1}{\text{min}}}

\makeatletter
\def\ps@pprintTitle{%
     \let\@oddhead\@empty
     \let\@evenhead\@empty
     \def\@oddfoot{\hfil}%
     \let\@evenfoot\@oddfoot}
\makeatother

\begin{document}

\begin{frontmatter}



\title{Multi-stage PDE-based image processing techniques for noisy MRI scans} 



\author[tud]{Ksenia~Slepova\corref{cor1}} 
\ead{k.slepova@tudelft.nl}
\cortext[cor1]{Corresponding author}

\author[jhsm]{Ivan~Etoku~Oiye}
\author[tud]{Martin~B.~van~Gijzen}

\affiliation[tud]{organization={Delft University of Technology, Delft Institute of Applied Mathematics},
            city={Delft},
            postcode={2628 CD}, 
            country={The Netherlands}}
\affiliation[jhsm]{organization={Accessible Magnetic Resonance Laboratory, The Russell H. Morgan Department of Radiology and Radiological Science, Johns Hopkins School of Medicine},
            city={Baltimore, Maryland},
            country={USA}}     

\begin{abstract}
Image denoising and image segmentation play essential roles in image processing. Partial differential equations (PDE)-based methods have proven to show reliable results when incorporated in both denoising and segmentation of images. In our work, we discuss a multi-stage PDE-based image processing approach. It relies upon the nonlinear diffusion for noise removal and clustering and region growing for segmentation. In the first stage of the approach, the raw image is computed from noisy measurement data. The second stage aims to filter out the noise using anisotropic diffusion. We couple these stages into one optimisation problem which allows us to incorporate a diffusion coefficient based on a presegmented image. The third stage performs the final segmentation of the image.  We demonstrate our approach on both images for which the ground truth is known and on MR measurements made by an experimental, inexpensive scanner.
\end{abstract}


\end{frontmatter}



\section{Introduction}
\label{S.1}

\subsection{Context of the research}

Image denoising and image segmentation play important roles in image processing. In many cases, images have to be reconstructed from noisy measurement data, such as magnetic resonance imaging (MRI) or computed tomography (CT) scans, resulting in noisy images that require further processing. 

This research was motivated by an earlier project that resulted in the development of a low-cost low-field MRI machine \cite{o2021vivo, o2019three, obungoloch2023site}. Images obtained from such scanners are often of poor quality and cannot be used by clinical doctors without additional processing. The aim of this work is, therefore, to develop a technique that would allow us to efficiently denoise and segment such images \cite{shan2022deflated}.

\subsection{Motivation}

Research on denoising and segmentation of noisy images often relies on PDE-based models \cite{weickert2001efficient, chen2000image, karasev2013interactive, kollem2019review, kollem2021, RIYA2021106}. However, nowadays it is to some extent overshadowed by artificial intelligence (AI) techniques \cite{ranjbarzadeh2021brain, liu2023deep, guo2019deep}. One of the disadvantages of employing AI-based methods is the lack of transparency in the decisions and actions the models make, whereas PDE-based methods are well explained and can be easily interpreted. This served as a motivation to utilize PDE-based approach in our research.

\subsection{Contributions of this paper}

This paper focuses on analytical, PDE-based techniques, extending the work in \cite{shan2022deflated} that explains computationally efficient noise reduction techniques for low-field MRI. In this paper, we discuss full image processing pipeline. Our approach is closely related to the ideas proposed by Wu \emph{et al.} in \cite{wu2021two}. Wu \emph{et al.~}suggest an approach that first uses an approximation model to filter the initial image, followed by segmentation through thresholding and k-means clustering \cite{macqueen1967some}.

Our processing pipeline comprises three major phases that are solved interconnectedly: image formation, image enhancement, and image segmentation. Firstly, we solve a minimization problem that includes both image filtering and image reconstruction terms. To ensure the efficiency of the method, we use the deflated preconditioned conjugate gradient method, as explained in \cite{shan2022deflated}. In the next stage, several enhancement techniques, such as background removal through region growing \cite{adams1994seeded}, morphological operations \cite{shapiro2001computer, gonzalez2006dip}, and histogram equalization \cite{gonzalez2006dip}, are applied, enriching the ideas from \cite{wu2021two} with a larger set of techniques. Finally, the Jenks natural breaks classification method \cite{jenks1967data} is employed to obtain the final segmentation of the initially noisy MRI scan. Moreover, we extend the approach by performing presegmentation prior to the filtering. We aim to investigate whether this multi-stage PDE-based image processing approach, with its tightly connected three phases, extended by adding a presegmentation step, yields more accurate segmentation results than the basic one.

\subsection{Structure of this paper}

In section \ref{S.2}, we will describe the filtering model that is used in our research. Next, in section \ref{S.3}, the basic multi-stage segmentation approach will be discussed, and the results of its performance will be presented. Section \ref{S.4} will cover the modified approach and will also include the results of segmentation with the use of it. Finally, some conclusions will be drawn in section \ref{S.5}.

\section{Image filtering technique}
\label{S.2}

Detailed derivations of the equations in this section can be found in \cite{shan2022deflated}.

\subsection{Diffusion PDE's for noise reduction}
\label{noise-reduction}
The Perona-Malik diffusion model is a filtering technique that allows to reduce image noise without diffusing the edges or other significant image details \cite{perona1990scale, weickert1998anisotropic}. It is given by the equations below:
\begin{align}
&\dd{u(\vecx,t)}{t}  = \nabla \cdot (c(\|\nabla u(\vecx,t)\|)\nabla u(\vecx,t)) ~~~ \text{in}~\Omega \times (0,T) \label{eq1}\\
&u(\vecx,0) = f(\vecx) ~~~\text{in}~\Omega \label{eq2}\\
&\dd{u(\vecx,t)}{\vecn} = 0 ~~~ \text{on}~\dg \times (0,T). \label{eq3}
\end{align}
Here, $\Omega$ is the picture domain, $T$ is the stopping time, $u(\vecx,t)$ is the image intensity, $f(\vecx)$ is the initial noisy image, and $c$ is a function satisfying following conditions:
\begin{itemize}
\item $c(\|\nabla u(\vecx,t)\|)$ is a non-negative function monotonically decreasing  in $\nabla u(\vecx,t)$,
\item $\lim\limits_{\|\nabla u(\vecx,t)\| \to  0} c(\|\nabla u(\vecx,t)\|) = M,$ where $M \in (0, \infty),$
\item $\lim\limits_{\|\nabla u(\vecx,t)\| \to  \infty} c(\|\nabla u(\vecx,t)\|) = 0$.
\end{itemize} 

The basic anisotropic diffusion PDE from Eq.(\ref{eq1}) can be modified in such a way that the updated equation maintains the fidelity to the original image,  effectively controlling how close the denoised image remains to the raw image. Moreover, in the modified model the stopping time does not have to be chosen any more and the diffusion termination at trivial solutions, such as a constant image, is avoided (see \cite{nordstrom1990biased}). Adding a fidelity term to the right-hand side which brings us to the following PDE:

\begin{align}
&\dd{u(\vecx,t)}{t}  = \nabla \cdot (c(\|\nabla u(\vecx,t)\|)\nabla u(\vecx,t)) + \eta(f(\vecx)- u(\vecx,t)) ~~~ \text{in}~\Omega \times (0,T) \label{eq8}
\end{align}
The fidelity term $\eta(u(\vecx, 0) - u(\vecx,t))$ in Eq.(\ref{eq8}) is providing a constraint that penalizes variation of the output image from the input one. This constraint is controlled by the fidelity parameter $\eta$.

\subsubsection{Diffusion coefficients}
\label{diff-coefs}

The original model suggests the following diffusion coefficients:
\begin{align}
c_1(\|\nabla u(\vecx,t)\|) &= \exp\left(-\left(\dfrac{\|\nabla u(\vecx,t)\|}{K}\right)^2\right), \label{eq4}\\
c_2(\|\nabla u(\vecx,t)\|) &= \dfrac{1}{1 + \left[\dfrac{\|\nabla u(\vecx,t)\|}{K}\right]^2}. \label{eq5}
\end{align}
Despite their ability to sharpen edges, these coefficients are unable to remove heavy-tailed noise and create the so-called ``staircase'' artefacts \cite{whitaker1993multi, you1994analysis}. Although this ``staircase'' effect is not useful for image enhancement, it might be beneficial to segmentation purposes.

The improved version of the original coefficient is a parameter-free total variation option:
\begin{align}
c_3(\|\nabla u(\vecx,t)\|) = \dfrac{1}{\|\nabla u(\vecx,t)\|}.  \label{eq7}
\end{align}

In our research,  we use the generalization of the elastic net regularization \cite{elasticnet}:
\begin{align}
c_{4}(\|\nabla u(\vecx,t)\|) = \dfrac{d_p}{2} \cdot \left(\dfrac{1}{\|\nabla u(\vecx,t)\|^2}\right)^{1-\frac{d_p}{2}} + K. \label{eq7b}
\end{align}
Here, $d_p$ is the diffusion power and $K$ is the diffusion parameter. 

Let us find the range of $d_p$ and $K$ such that the suggested choice of the diffusion coefficient  is non-negative, monotonically decreasing in $\nabla u(\vecx,t)$, has a finite non-zero limit when $\nabla u(\vecx,t) \to 0$, and has a zero limit when $\nabla u(\vecx,t)$ approaches infinity. For this, we rewrite the diffusion coefficient in the following way:
\begin{align}
c_{4}(\|\nabla u(\vecx,t)\|) = \dfrac{d_p}{2} \cdot \left(\|\nabla u(\vecx,t)\|\right)^{2 \cdot (\frac{d_p}{2}-1)} + K = \dfrac{d_p}{2} \cdot \left(\|\nabla u(\vecx,t)\|\right)^{d_p-2} + K \label{eq7c}. 
\end{align}
We now consider only the first term of the diffusion coefficient, as the second one is independent of $\|\nabla u\|$. From the expression in (\ref{eq7c}), we conclude that the first term is monotonically decreasing as $\|\nabla u\|$ increases when $d_p < 2$.

Now, let us analyse the behaviour of the diffusion coefficient as $\|\nabla u\| \to 0$. With $d_p < 2$, the power $d_p - 2 < 0$, therefore, the first term of $c_4$ blows up. To prevent this, we add a small regularization term $\eps$, which brings us to the following form of the diffusion coefficient:
\begin{align}
c_{4}(\|\nabla u(\vecx,t)\|) = \dfrac{d_p}{2} \cdot \left(\dfrac{1}{\|\nabla u(\vecx,t)\|^2 + \eps}\right)^{1-\frac{d_p}{2}} + K. \label{eq7d}
\end{align}
In our experiments, we take $\eps = 10^{-4}$. Finally, when the gradient approaches infinity, we get the following limit:
\begin{align}
\lim\limits_{\|\nabla u(\vecx,t)\| \to \infty} c_{4}(\|\nabla u(\vecx,t)\|) = K. \label{eq7e}
\end{align}
We aim to make this limit as close to $0$ as possible but keep it non-negative. As we scale all the images between $0$ and $1$, we will only consider positive $K$ of order $10^{-1}$ and below that.

\subsection{Picard iteration}
\label{picard-iter}
In \cite{vese2015variational}, the authors explain that, under certain conditions, the solution $u(\vecx, t)$ of the time-dependent problem should approximate a minimizer $u(\vecx)$ of the model (\ref{eq8b}) as time increases. Therefore, we will find the stationary solution $u(\vecx)$ of (\ref{eq8}). It should satisfy the following equation:
\begin{align}
&0  = \nabla \cdot \left(c(\|\nabla u(\vecx)\|)\nabla u(\vecx) \right) + \eta(f(\vecx) - u(\vecx)) ~~~ \text{in}~\Omega \label{eq8b}
\end{align}

To solve the system, we firstly discretize Eq.(\ref{eq8b}) in space using the standard finite difference method (FDM) (see, for instance, \cite{perona1990scale, catte1992image}). We obtain the following semi-discrete system:
\begin{align}
0 = C(\vecu) \vecu + \eta ( \vecf - \vecu), \label{eq8d}
\end{align}
where $\vecu$ and $\vecf$ contain pixel values of the denoised image and of the original noisy image, respectively.

Rewriting Eq.(\ref{eq8d}) brings us to the following form:
\begin{align}
\left(I - \frac{1}{\eta}C(\vecu)\right)\vecu = \vecf. \label{eq8e}
\end{align}

The non-linear system in Eq.(\ref{eq8e}) is then solved using the lagged diffusion Picard iteration as suggested in \cite{vogel1996total}:
\begin{align}
\left(I - \frac{1}{\eta}C\left(\vecu^{n}\right)\right)\vecu^{n+1} = \vecf. \label{system-final}
\end{align}

Denoting $A(\vecu) := \left(I - \frac{1}{\eta}C(\vecu)\right), ~ \vecb := \vecf$, we obtain the final linear system that needs to be solved in every Picard iteration:
\begin{align}
A\left(\vecu^n\right)\vecu^{n+1} = \vecb. \label{system-final-picard}
\end{align}

Matrix $A$ in Eq.(\ref{system-final-picard}) is symmetric and positive definite. However, discontinuities in the diffusion coefficient might occur, resulting in $A$ being ill-conditioned. System (\ref{system-final-picard}) can be solved using the preconditioned conjugate gradient method:
\begin{align}
M^{-1}A \vecu = M^{-1}\vecb,
\end{align}
where $M$ is chosen in such a way that it resembles $A$. Further improvements in the speed of the convergence can be performed by applying deflated preconditioner (DPCG). That allows to map isolated eigenvalues to zero, effectively removing them from the system. Preconditioners and deflation techniques for system \ref{system-final-picard} are discussed in detail in \cite{shan2022deflated} . In our model we use subdomain deflation as described in \cite{nicolaides1987deflation}.

\section{Basic multi-stage image processing approach}
\label{S.3}

\subsection{Image reconstruction}
MRI raw data is collected in the frequency domain, therefore, it needs to be converted into the spatial domain for further processing. This process is called image reconstruction, or image formation, and it is performed by applying an inverse Fourier transform \cite{gallagher2008introduction}. Thus, the relation between the initial noisy frequency-domain data $f$ and the spatial image $u_0$ can be written as follows:
\begin{align}
\mathcal{A}u_0 =  f,
\end{align}
where $\mathcal{A}$ represents the Fourier transform operator.

\subsection{Denoising}

Consider the following minimization problem:
\begin{align}
\argmin{u} ~ \left[ \frac{\eta}{2} \int\limits_{\Omega} (u_0- u)^2 ~dx  +  R(\|\nabla u\|^2) \right]. \label{integral-model}
\end{align}
Here,  $\eta$ is a fidelity parameter, $u$ is an optimal piecewise smooth approximation of the initial image $u_0$, $\Omega \subset \real^2$ is the image domain. Moreover, $\Omega$ is supposed to be a bounded open set with Lipschitz boundary and $u_0$ is continuous inside the image domain. Note that $u_0$ is the result of the inverse Fourier transform of an MRI signal. The second term is the regularization term. It should filter out noise while preserving the edges of the objects.

The energy functional given by the expression in (\ref{integral-model}) is minimized by the solution of the corresponding Euler-Lagrange equations:
\begin{align}
\nabla \cdot \left(R'(\|\nabla u\|)\nabla u\right) + \eta(u_0 - u) = 0. \label{euler-lagrange}
\end{align}
One can notice that it coincides with the earlier derived Eq.(\ref{eq8b}) because $R'(\cdot)$ is equal to the diffusion coefficient $c(\cdot)$. A detailed discussion of possible choices for $c(\cdot)$ can be found in section \ref{diff-coefs}. In Table (\ref{regularise-matrix}), we present only total variation (TV) and generalized elastic net regularizers. The latter is used in our experiments.

\renewcommand{\arraystretch}{2.0}
\begin{table}[H]
\captionsetup{justification=centering}
\centering
\caption{Choice of regularizers. Note that $s = \|\nabla u \|^2$.}
\begin{tabular}{l||c c}
\hline 
Name & $R(s)$  & $R'(s) = c(s)$ \\ 
\hline \hline
Total Variation & $2 \sqrt{s}$ & $\frac{1}{\sqrt{s}}$ \\
Generalized Elastic Net & $s^{\frac{d_p}{2}} + Ks$ & $~~~\frac{d_p}{2} \cdot \left(\dfrac{1}{s}\right)^{1-\frac{d_p}{2}} + K$\\
\hline
\end{tabular}
\label{regularise-matrix}
\end{table}

\subsection{Minimization problem}

Image reconstruction and filtering can be combined in one optimization problem, since $R(\cdot)$ can serve as a regularizer for $f$ as well:
\begin{align}
\min\limits_{v, u} ~ \left[ h(v) + g(u) \right] ~~~ \text{s.t.} ~ u - v = 0. \label{discrete-model}
\end{align}

Here,  
\begin{enumerate}
\item Image reconstruction is given by:  \begin{equation}
h(v) = \frac{\eta}{2} \|\mathcal{A}v - f\|^2_2; ~~~\min\limits_{v} \left[  h(v) \right] ; \label{sub1}
\end{equation} $\mathcal{A}$ is the Fourier transform operator \cite{de2022image}; expressions in (\ref{sub1}) come from the following derivation:
\begin{equation*}
\min\limits_{v} \|v - u_0\|^2_2 \stackrel{\text{orthogonal}~\mathcal{A}}{=} \min\limits_{v} \|\mathcal{A}v - \mathcal{A}u_0\|^2_2  = \min\limits_{v} \|\mathcal{A}v - f\|^2_2;
\end{equation*}
\item Nonlinear diffusion filtering is represented by:  \begin{equation}
g(u) = R(\|\nabla u\|^2); ~~~\min\limits_{u}  \left[  g(u) \right] ;
\end{equation}
\end{enumerate}
This minimization problem is then solved using the alternating direction method of multipliers (ADMM) \cite{boyd2011distributed}.

\subsection{Segmentation}
\label{kmeans}

The main idea of the $k$-means clustering method is to give a partition of the initial set into $k$ disjoint clusters (with $k$ being a predefined parameter) \cite{macqueen1967some}. In our research, we use the Jenks natural breaks classification method, which is a computationally efficient variation of the k-means clustering method but applied to univariate data \cite{jenks1967data}. This is an iterative process that consists of the following steps:
\begin{enumerate}
\item Choose the maximum number of clusters $n_{\text{cl.}}$;
\item Split the set into $k$ initial clusters in some way which can be arbitrary;
\item Calculate the sum of squared deviations from the class means (SDCM);\label{step1}
\item Regroup the data into new clusters, possibly by moving elements from one cluster to a different one; \label{step2}
\item Compute the new sum of deviations per cluster. \label{step3}
\end{enumerate}
Steps \ref{step1}, \ref{step2}, and \ref{step3} are repeated until a certain tolerance $\tau$ is achieved or SDCM becomes constant. 

Let us specify the way in which the set is divided into clusters and which stopping criterion is used in our model. Firstly, all the pixel values are sorted, and the number of clusters $k$ is set to $1$. After that, SDCM is computed. Moreover, the sum of squared deviations from the mean of the whole dataset (SDAM) and the goodness of variance fit (GVF) are calculated. The number of clusters is iteratively increased, and the elements are moved from one class to another while 
\begin{align*}
k < n_{\text{cl.}} ~ \text{and} ~ \text{GVF} = \dfrac{\text{SDAM} - \text{SDCM}}{\text{SDAM}} < \tau.
\end{align*}

\subsection{Multi-stage approach}
\label{basic_multistage}

Finally, let us describe the basic approach for filtering and segmentation of noisy MRI scans. The method has several stages. In the first stage, the minimization problem (\ref{discrete-model}) is solved; thus, an image is formed, and an approximation of $u$ is obtained. In order to perform that, we use a single step of ADMM.  

In the second stage, we perform splitting of the solution $u$ into clusters. Based on the idea of Wu \emph{et al.} \cite{wu2021two}, we then use the following basic multi-stage approach:
\begin{enumerate}
\item First stage: 
\begin{itemize}
\item Form the initial image by reconstructing it from the existing data;
\item Enhance the image using anisotropic diffusion filtering with elastic net as a diffusion coefficient; 
\end{itemize}
\item Second stage: apply (optionally) background removal through region growing \cite{adams1994seeded}, morphological closing \cite{shapiro2001computer, gonzalez2006dip}, global or adaptive histogram equalization to the filtered image \cite{gonzalez2006dip},  and then segment it using the Jenks natural breaks classification method (see section \ref{kmeans}).
\end{enumerate}

\subsection{Results}

\subsubsection{Evaluation criteria}
\label{metrics}

In numerical experiments, we will compare the results with the ground truth (GT) using several quantitative indicators. For the images that do not have GT, we will evaluate the results visually. The evaluation criteria that we will be using are the Jaccard score (JS), Dice similarity coefficient (DSC), and segmentation accuracy (SA) \cite{muller2022towards}. Let $A$ denote the region segmented by our approach and $B$ denote corresponding region in the GT segmentation. The aforementioned evaluation criteria can then be defined as follows:

\begin{enumerate}
\item Jaccard score (JS), or Intersection-over-Union (IoU):
\begin{equation*}
JS(A,B) = \dfrac{|A \cap B|}{|A \cup B|}.
\end{equation*}
\item Dice similarity coefficient (DSC):
\begin{equation*}
DSC(A,B)= \dfrac{2 \cdot |A \cap B|}{|A|+|B|}.
\end{equation*}
\item Segmentation accuracy (SA):
\begin{equation*}
SA(A,B) = \dfrac{|A \cap B| + |A^c \cap B^c|}{|A|+|A^c|},
\end{equation*}
where $A^c$ and $B^c$ are the complements of $A$ and $B$ respectively.
\end{enumerate}

\subsubsection{Testing dataset}
\label{testing-set}

\begin{figure}[H]
    \begin{subfigure}{.32 \textwidth}
        \centering
        \includegraphics[width=2.5cm]{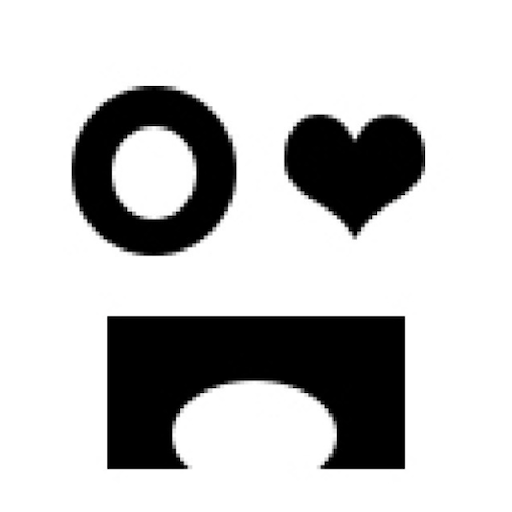}
        \caption{\centering Given image}
        \label{1a}
    \end{subfigure}
    \begin{subfigure}{.32 \textwidth}
        \centering
        \includegraphics[width=2.5cm]{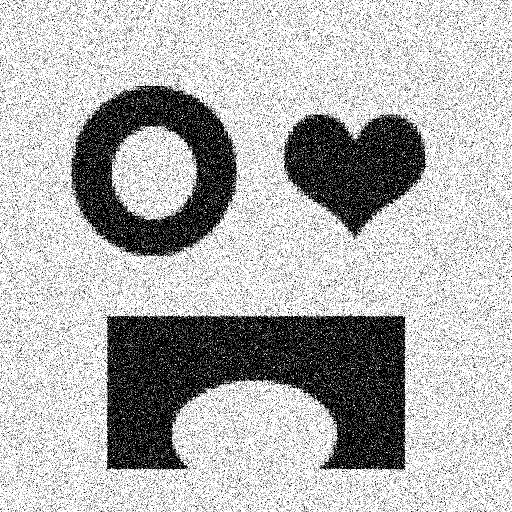}
        \caption{\centering $\sigma^2 = 0.1$}
		\label{1b}
    \end{subfigure}
    \begin{subfigure}{.32 \textwidth}
        \centering
        \includegraphics[width=2.5cm]{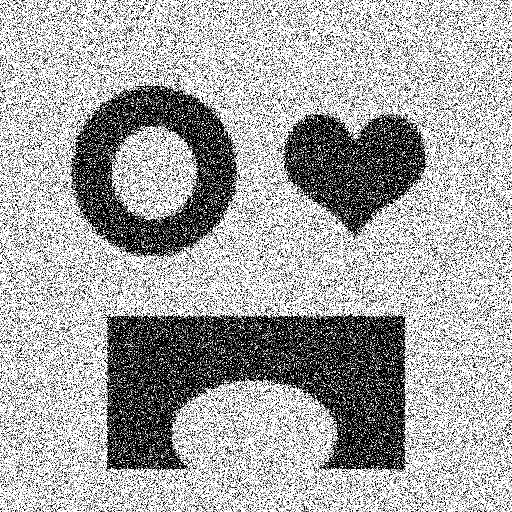}
        \caption{\centering $\sigma^2 = 0.3$}
		\label{1c}
    \end{subfigure}\\[0.1cm]
    \begin{subfigure}{.32 \textwidth}
        \centering
        \includegraphics[width=2.5cm]{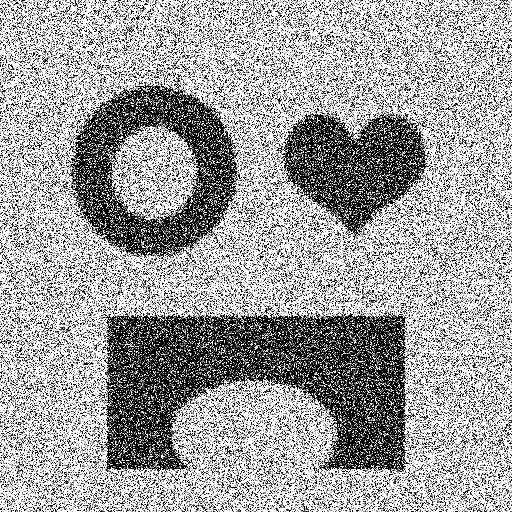}
        \caption{\centering $\sigma^2 = 0.5$}
		\label{1d}
    \end{subfigure}
    \begin{subfigure}{.32 \textwidth}
        \centering
        \includegraphics[width=2.5cm]{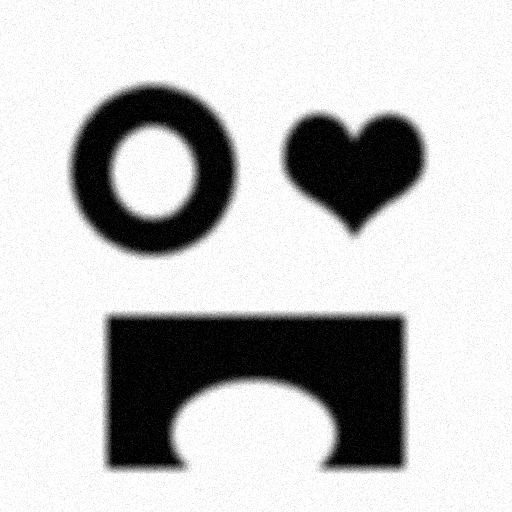}
        \caption{\centering $\mathcal{B}_g$}
		\label{1e}
    \end{subfigure}
    \begin{subfigure}{.32 \textwidth}
        \centering
        \includegraphics[width=2.5cm]{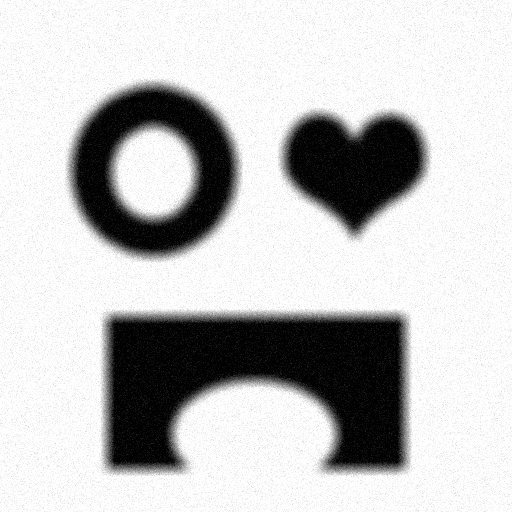}
        \caption{\centering $\mathcal{B}_a$}
		\label{1f}
    \end{subfigure}\\[0.1cm]
    \begin{subfigure}{.32 \textwidth}
        \centering
        \includegraphics[width=2.5cm]{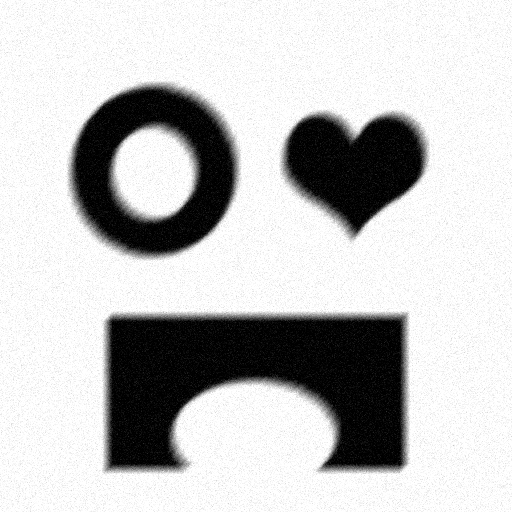}
        \caption{\centering $\mathcal{B}_m$}
		\label{1g}
    \end{subfigure}
    \begin{subfigure}{.32 \textwidth}
        \centering
        \includegraphics[width=2.5cm]{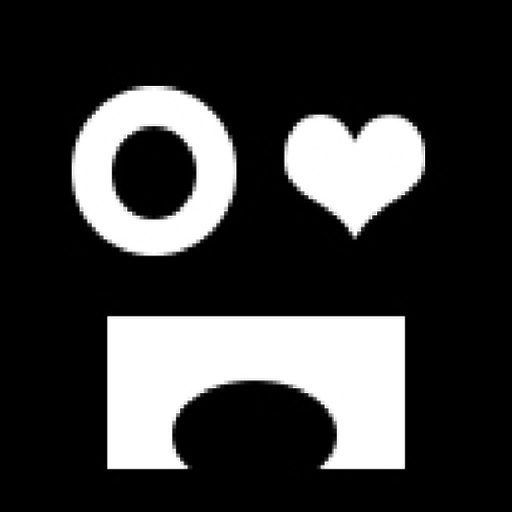}
        \caption{\centering Ground truth}
		\label{1h}
    \end{subfigure}
    \begin{subfigure}{.32 \textwidth}
        \centering
        \includegraphics[width=2.5cm]{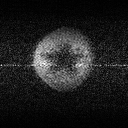}
        \caption{\centering Papaya scan}
		\label{1i}
    \end{subfigure}
    \caption{Raw images and scans used for testing.}
    \label{fig1}
\end{figure}

In this paper, we test our approach using two types of data: a synthetically generated image and an MRI scan from a low-field scanner. The synthetic image is presented in Fig. \ref{fig1}(a) and has a resolution of $512\times512$ pixels. The same image was used in the paper \cite{wu2021two}, and we will perform a  comparison of our results with the results of Wu \textit{et al.}. Gaussian noise with a mean of $0$ and variances of $0.1,~0.3,~0.5$ is added to the image, as shown in Fig. \ref{fig1}(b)-\ref{fig1}(d), respectively. Additionally, we test our approach on images with different types of blurring (see Fig. \ref{fig1}(e)-\ref{fig1}(g)). The blurring is applied to the image in MATLAB using the following  functions: $\mathcal{B}_g$=fspecial(`gaussian', $[12,12]$, $12$); $\mathcal{B}_a$=fspecial(`average', $15$); $\mathcal{B}_m$=fspecial(`motion', $15$, $45$); $\mathcal{B}_g, \mathcal{B}_a$ and $\mathcal{B}_m$ stand for Gaussian blur, average blur and motion blur, respectively. The ground truth segmentation for the image is shown in Fig. \ref{fig1}(h). 

Finally, an MRI scan of a papaya, acquired from a low-cost scanner located at Mbarara University of Technology \cite{obungoloch2023site}, is presented in Fig. \ref{fig1}(i). The scanner was still under development when the measurements were obtained; therefore, the results are preliminary. No ground truth segmentation is available for this scan.

\subsubsection{Numerical results}

\begin{figure}[H]
    \begin{subfigure}{.32 \textwidth}
        \centering
        \includegraphics[width=2.5cm]{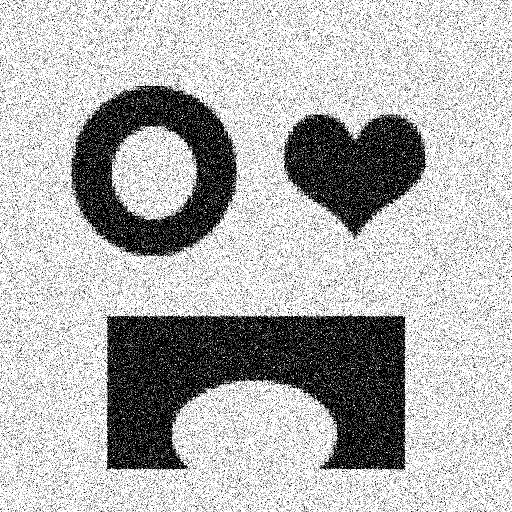}
        \caption{\centering $\sigma^2 = 0.1$}
    \end{subfigure}
    \begin{subfigure}{.32 \textwidth}
        \centering
        \includegraphics[width=2.5cm]{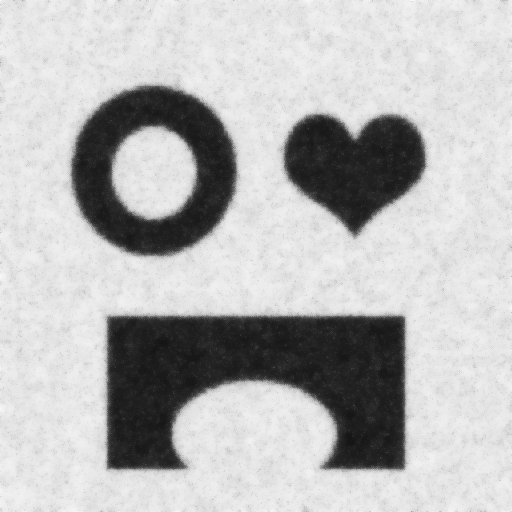}
        \caption{\centering $\sigma^2 = 0.1$}
    \end{subfigure}
    \begin{subfigure}{.32 \textwidth}
        \centering
        \includegraphics[width=2.5cm]{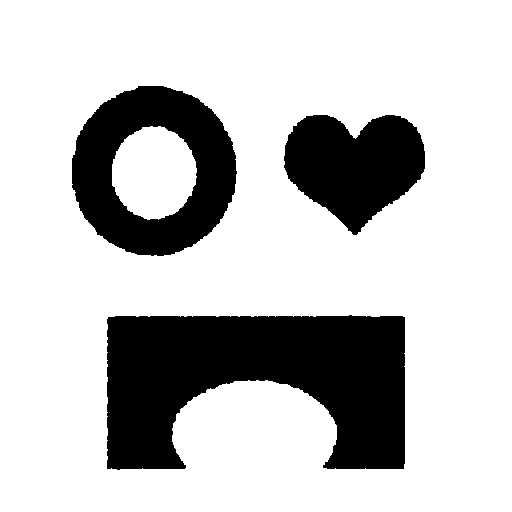}
        \caption{\centering $\sigma^2 = 0.1$}
    \end{subfigure}\\[0.1cm]
     \begin{subfigure}{.32 \textwidth}
        \centering
        \includegraphics[width=2.5cm]{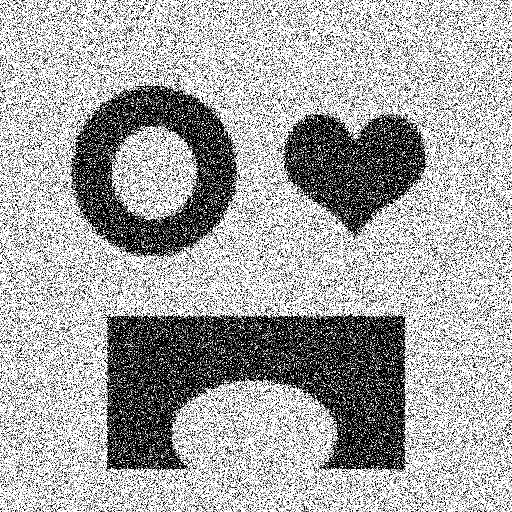}
        \caption{\centering $\sigma^2 = 0.3$}
    \end{subfigure}
    \begin{subfigure}{.32 \textwidth}
        \centering
        \includegraphics[width=2.5cm]{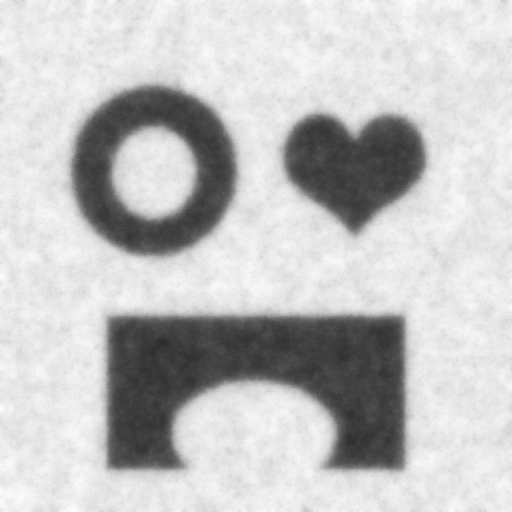}
        \caption{\centering $\sigma^2 = 0.3$}
    \end{subfigure}
    \begin{subfigure}{.32 \textwidth}
        \centering
        \includegraphics[width=2.5cm]{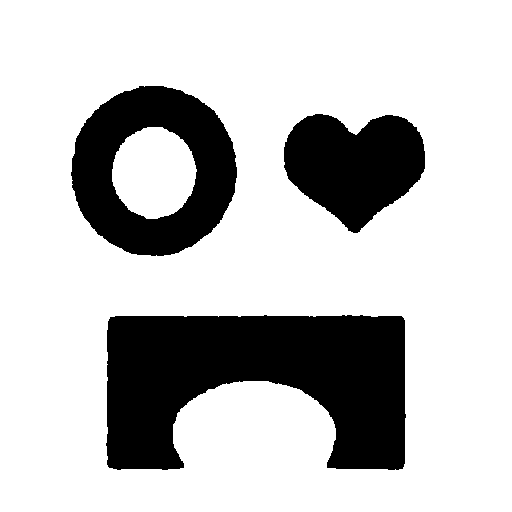}
        \caption{\centering $\sigma^2 = 0.3$}
    \end{subfigure}\\[0.1cm]
    \begin{subfigure}{.32 \textwidth}
        \centering
        \includegraphics[width=2.5cm]{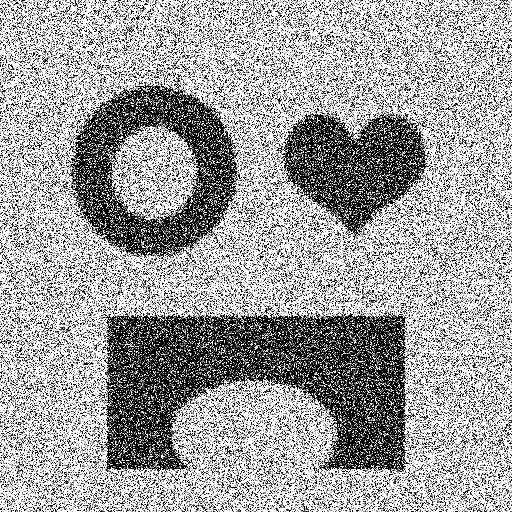}
        \caption{\centering $\sigma^2 = 0.5$}
    \end{subfigure}
    \begin{subfigure}{.32 \textwidth}
        \centering
        \includegraphics[width=2.5cm]{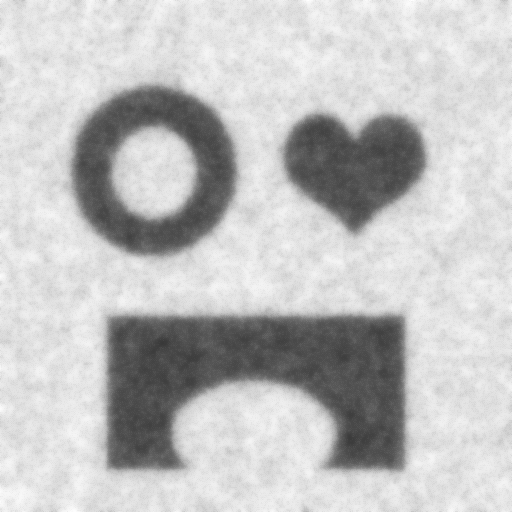}
        \caption{\centering $\sigma^2 = 0.5$}
    \end{subfigure}
    \begin{subfigure}{.32 \textwidth}
        \centering
        \includegraphics[width=2.5cm]{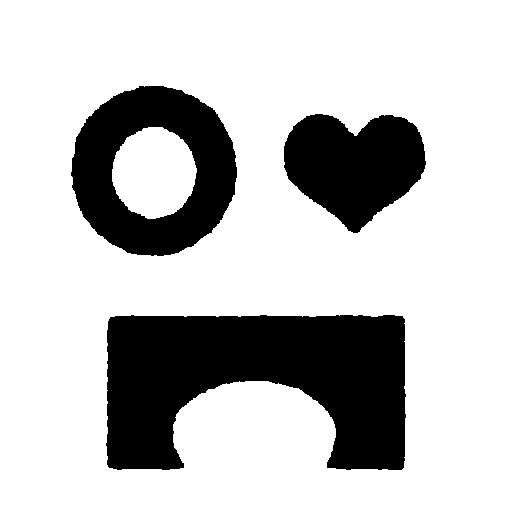}
        \caption{\centering $\sigma^2 = 0.5$}
    \end{subfigure}
    \caption{Results for noisy images using basic approach. \textit{Left column}: raw images, \textit{middle column}: filtered images, \textit{right column}: segmented images.}
    \label{fig2}
\end{figure}

\begin{figure}[H]
    \begin{subfigure}{.32 \textwidth}
        \centering
        \includegraphics[width=2.5cm]{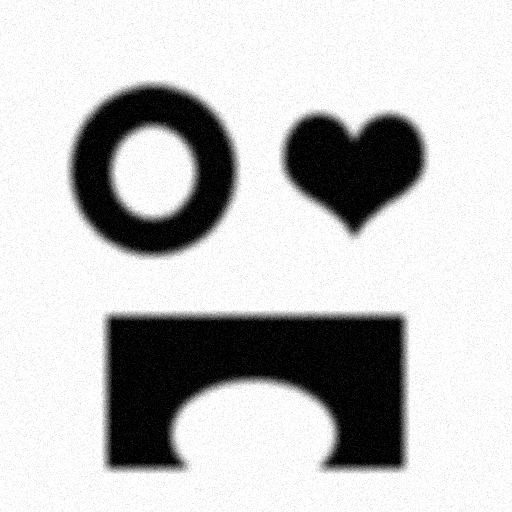}
        \caption{\centering $\mathcal{B}_g$}
    \end{subfigure}
    \begin{subfigure}{.32 \textwidth}
        \centering
        \includegraphics[width=2.5cm]{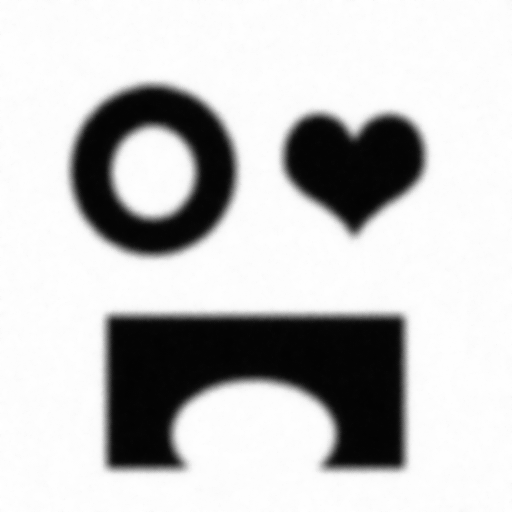}
        \caption{\centering $\mathcal{B}_g$}
    \end{subfigure}
    \begin{subfigure}{.32 \textwidth}
        \centering
        \includegraphics[width=2.5cm]{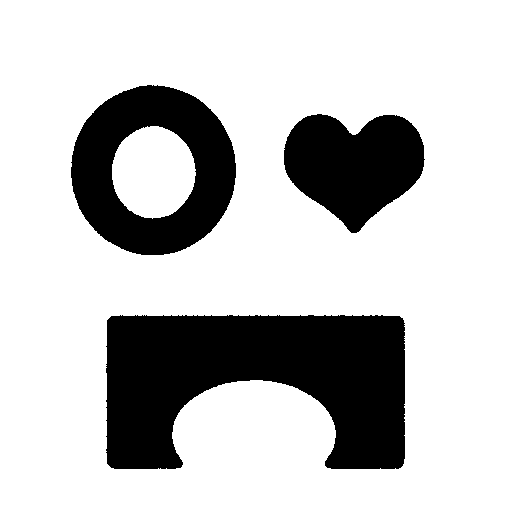}
        \caption{\centering $\mathcal{B}_g$}
    \end{subfigure}\\[0.1cm]
     \begin{subfigure}{.32 \textwidth}
        \centering
        \includegraphics[width=2.5cm]{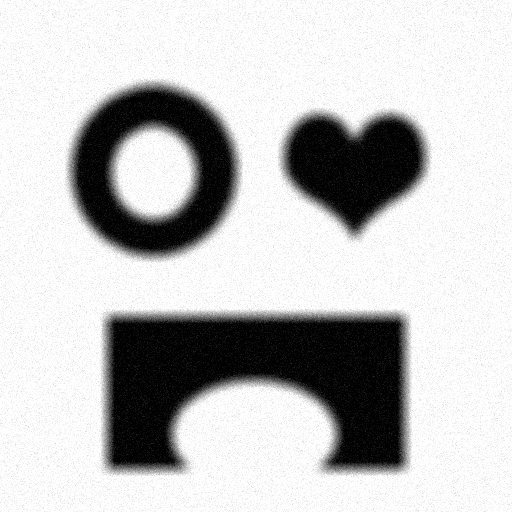}
        \caption{\centering $\mathcal{B}_a$}
    \end{subfigure}
    \begin{subfigure}{.32 \textwidth}
        \centering
        \includegraphics[width=2.5cm]{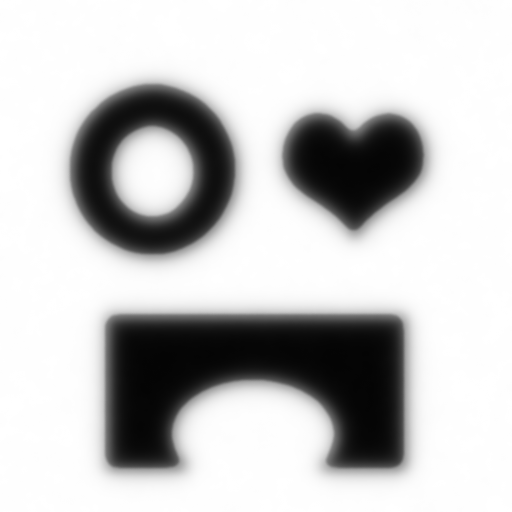}
        \caption{\centering $\mathcal{B}_a$}
    \end{subfigure}
    \begin{subfigure}{.32 \textwidth}
        \centering
        \includegraphics[width=2.5cm]{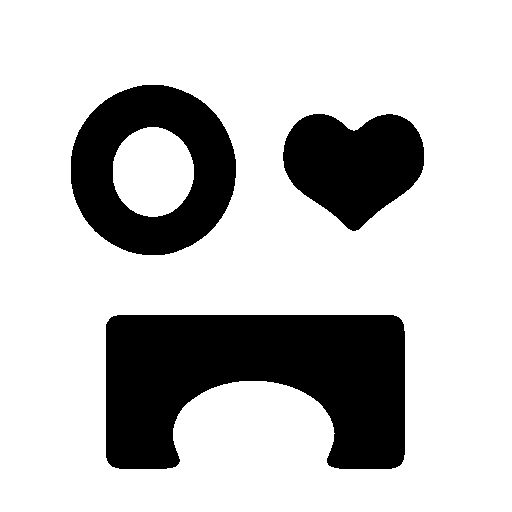}
        \caption{\centering $\mathcal{B}_a$}
    \end{subfigure}\\[0.1cm]
    \begin{subfigure}{.32 \textwidth}
        \centering
        \includegraphics[width=2.5cm]{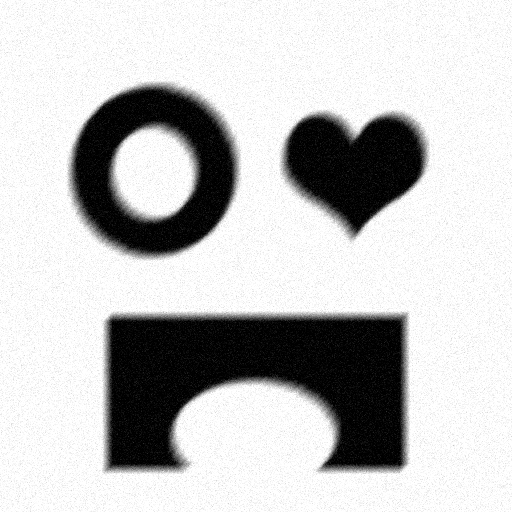}
        \caption{\centering $\mathcal{B}_m$}
    \end{subfigure}
    \begin{subfigure}{.32 \textwidth}
        \centering
        \includegraphics[width=2.5cm]{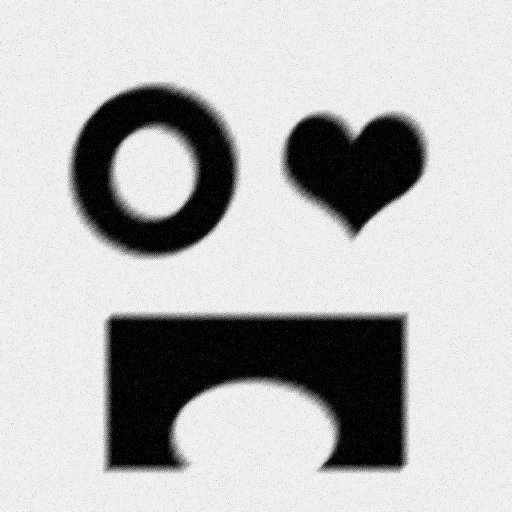}
        \caption{\centering $\mathcal{B}_m$}
    \end{subfigure}
    \begin{subfigure}{.32 \textwidth}
        \centering
        \includegraphics[width=2.5cm]{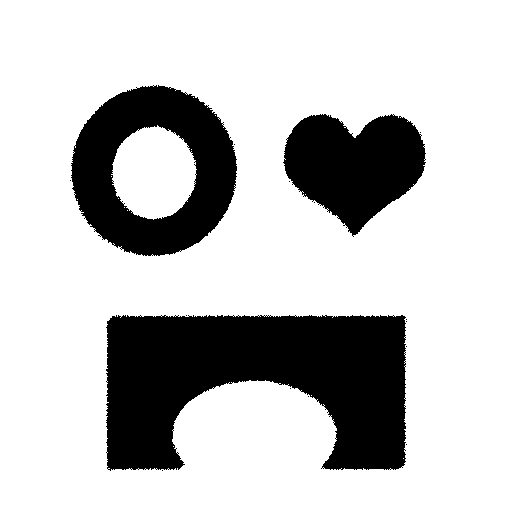}
        \caption{\centering $\mathcal{B}_m$}
    \end{subfigure}
    \caption{Results for blurry images using basic approach. \textit{Left column}: raw images, \textit{middle column}: filtered images, \textit{right column}: segmented images.}
    \label{fig3}
\end{figure}

We compare the performance of our approach  with that of the two-step image segmentation approach (TSIS) in \cite{wu2021two}. To do so, we compute the values of various quantitative metrics discussed in Section \ref{metrics}. The results obtained using the basic approach are shown in Fig. \ref{fig2} and Fig. \ref{fig3}.  In Fig. \ref{fig2}(c),(f),(i) and Fig. \ref{fig3}(c),(f),(i), white and black colours denote the two clusters (or the background and the foreground). Table \ref{noise1-compare} and Table \ref{blur1-compare} present a comparison of the evaluation metrics. Both visual and quantitative analyses indicate that the basic approach performs robustly across images with varying noise levels and types of blurring. Furthermore, quantitative analysis shows that our approach outperforms TSIS in some cases, particularly when $\sigma^2 = 0.5$.

\begin{figure}[H]
    \begin{subfigure}{.32 \textwidth}
        \centering
        \includegraphics[width=2.5cm]{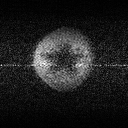}
    \end{subfigure}
    \begin{subfigure}{.32 \textwidth}
        \centering
        \includegraphics[width=2.5cm]{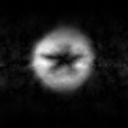}
    \end{subfigure}
    \begin{subfigure}{.32 \textwidth}
        \centering
        \includegraphics[width=2.5cm]{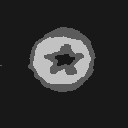}
    \end{subfigure}
    \caption{Results for the papaya MRI scan using basic approach. \textit{Left to right}: raw image, filtered image, segmented image.}
    \label{fig4}
\end{figure}

The results for the MRI scan of a papaya are shown in Fig. \ref{fig4}. Black represents the background and cavities within the object, dark gray corresponds to the peel and the star-shaped core, while the rest is indicated in light gray. Due to significant noise and the presence of a zipper artefact in the scan, we lowered the fidelity parameter $\eta$ compared to previous examples (see \ref{app1}). Analysis of the final segmentation demonstrates that the basic approach effectively handles heavy noise while preserving object edges, resulting in reliable outcomes.

\renewcommand{\arraystretch}{1.5}
\begin{table}[H]
\captionsetup{justification=centering}
\caption{Basic approach, images with different levels of Gaussian noise added (in $\%$).}
\centering
\begin{tabular}{ m{2cm}  m{1.2cm} m{1.2cm}  m{1.2cm} m{1.2cm}  m{1.2cm} m{1.2cm}}
\hline
& \multicolumn{2}{c}{JS} & \multicolumn{2}{c}{DSC}&  \multicolumn{2}{c}{SA} \\
\hline
Image & TSIS & Ours & TSIS & Ours & TSIS & Ours\\
\hline
$\sigma^2 = 0.1$ & \textbf{98.9} & 98.52 & \textbf{99.4} & 99.25 & \textbf{99.7} & 99.66\\
$\sigma^2 = 0.3$ & 97.6 & \textbf{97.80} & 98.9 & \textbf{98.89} & 99.4 & \textbf{99.49}\\
$\sigma^2 = 0.5$ & 96.53 & \textbf{97.51} & 98.2 & \textbf{98.74} & 99.2 & \textbf{99.42}\\
\hline
\end{tabular}
\label{noise1-compare}
\end{table}

\renewcommand{\arraystretch}{1.5}
\begin{table}[H]
\captionsetup{justification=centering}
\caption{Basic approach, images with different types of blurring (in $\%$).}
\centering
\begin{tabular}{ m{2cm}  m{1.2cm} m{1.2cm}  m{1.2cm} m{1.2cm}  m{1.2cm} m{1.2cm}}
\hline
& \multicolumn{2}{c}{JS} & \multicolumn{2}{c}{DSC}&  \multicolumn{2}{c}{SA} \\
\hline
Image & TSIS & Ours & TSIS & Ours & TSIS & Ours\\
\hline
$\mathcal{B}_g$ & \textbf{98.9} & 97.69 & \textbf{98.9} & 98.43 & \textbf{99.5} & 99.46\\
$\mathcal{B}_a$ & \textbf{98.3} & 97.74 & 97.7 & \textbf{98.86} & 98.9 & \textbf{99.47}\\
$\mathcal{B}_m$ & \textbf{99.4} & 97.48 & \textbf{99.5} & 98.73 & \textbf{99.9} & 99.41\\
\hline
\end{tabular}
\label{blur1-compare}
\end{table}

\section{Modified multi-stage image processing approach}
\label{S.4}
In the basic approach discussed in section~\ref{S.3}, the diffusion coefficient is computed directly from the raw image. In this section, we introduce a modified approach. Our observations indicate that computing the diffusion coefficient based on a presegmented image, combined with Gaussian smoothing (see \ref{app2}), produces sharper transitions in the diffusion coefficients. Using this observation, we incorporate a presegmentation step into the model. This modification extends the work of Wu \emph{et al.}~\cite{wu2021two} by integrating segmentation into the denoising process and iteratively performing segmentation and filtering iteratively at this stage. The modified approach thus consists of the following stages:

\begin{enumerate}
\item First stage: 
\begin{itemize}
\item Form the initial image by reconstructing it from the existing data;
\item Perform segmentation of the original image, apply Gaussian smoothing \cite{gonzalez2006dip} to the presegmented image and compute the diffusion coefficient based on that image;
\item Enhance the original image using anisotropic diffusion filtering and the diffusion coefficient computed earlier; 
\end{itemize}
\item Second stage: apply (optionally) background removal, morphological closing, global or adaptive histogram equalization to the filtered image, and then segment it using the Jenks natural breaks classification method.
\end{enumerate}

\subsection{Numerical results}

\begin{figure}[h]
    \begin{subfigure}{.32 \textwidth}
        \centering
        \includegraphics[width=2.5cm]{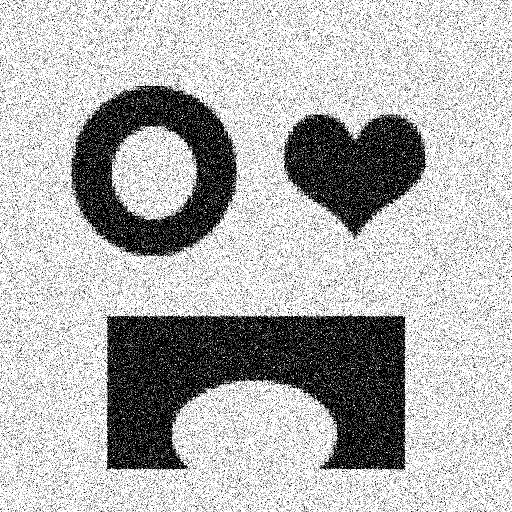}
        \caption{\centering $\sigma^2 = 0.1$}
    \end{subfigure}
    \begin{subfigure}{.32 \textwidth}
        \centering
        \includegraphics[width=2.5cm]{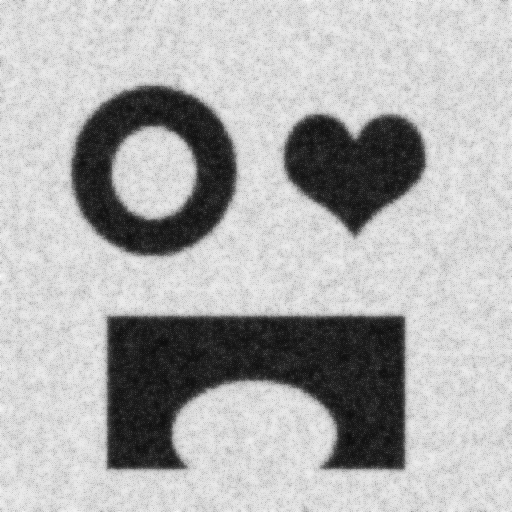}
        \caption{\centering $\sigma^2 = 0.1$}
    \end{subfigure}
    \begin{subfigure}{.32 \textwidth}
        \centering
        \includegraphics[width=2.5cm]{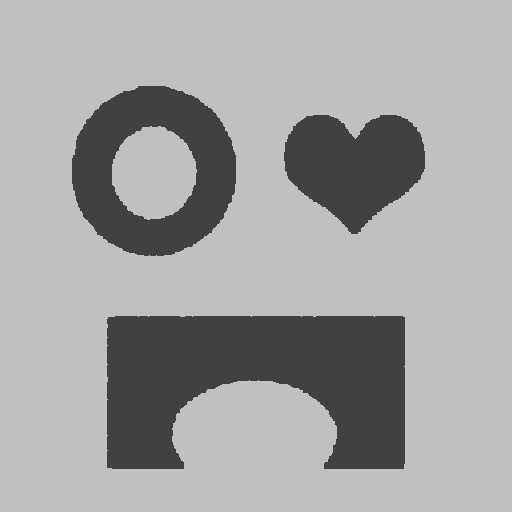}
        \caption{\centering $\sigma^2 = 0.1$}
    \end{subfigure}\\[0.1cm]
     \begin{subfigure}{.32 \textwidth}
        \centering
        \includegraphics[width=2.5cm]{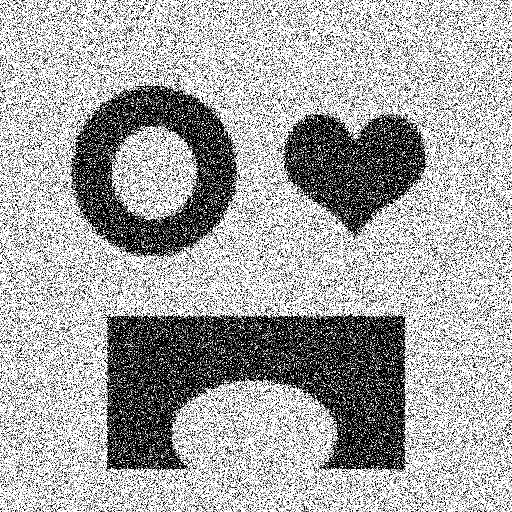}
        \caption{\centering $\sigma^2 = 0.3$}
    \end{subfigure}
    \begin{subfigure}{.32 \textwidth}
        \centering
        \includegraphics[width=2.5cm]{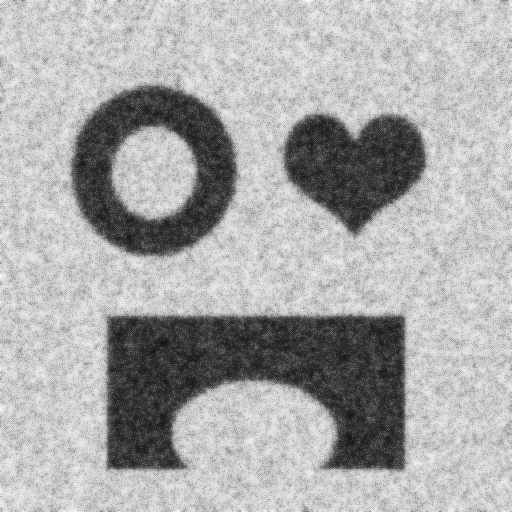}
        \caption{\centering $\sigma^2 = 0.3$}
    \end{subfigure}
    \begin{subfigure}{.32 \textwidth}
        \centering
        \includegraphics[width=2.5cm]{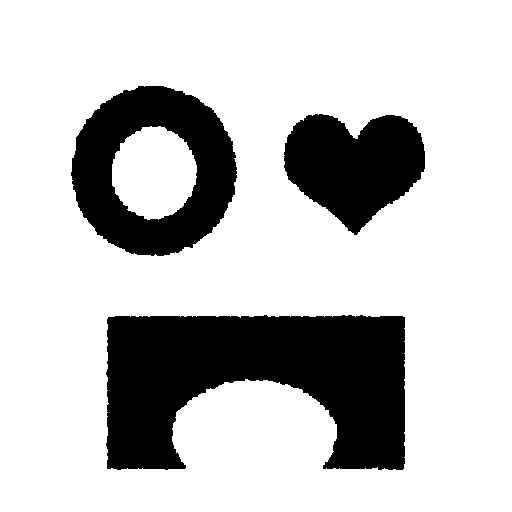}
        \caption{\centering $\sigma^2 = 0.3$}
    \end{subfigure}\\[0.1cm]
    \begin{subfigure}{.32 \textwidth}
        \centering
        \includegraphics[width=2.5cm]{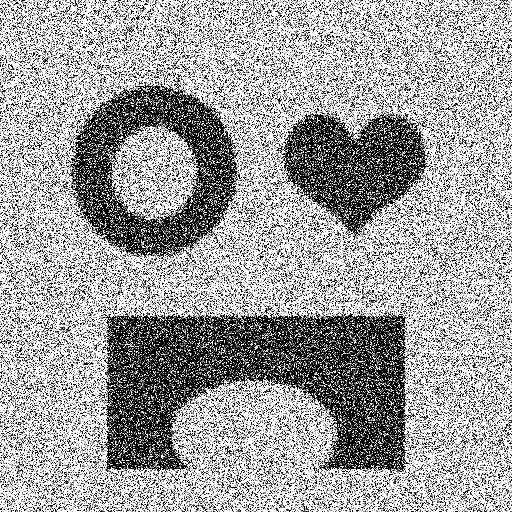}
        \caption{\centering $\sigma^2 = 0.5$}
    \end{subfigure}
    \begin{subfigure}{.32 \textwidth}
        \centering
        \includegraphics[width=2.5cm]{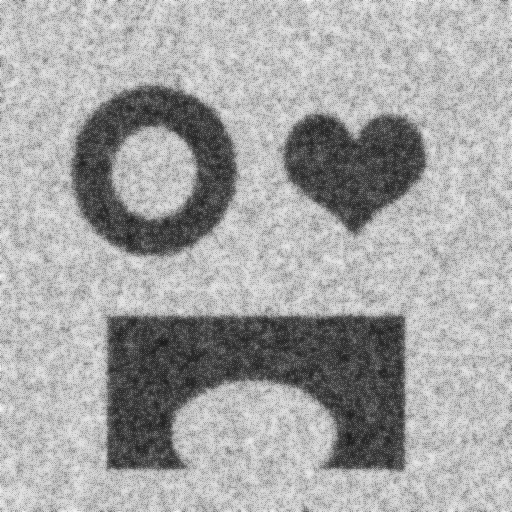}
        \caption{\centering $\sigma^2 = 0.5$}
    \end{subfigure}
    \begin{subfigure}{.32 \textwidth}
        \centering
        \includegraphics[width=2.5cm]{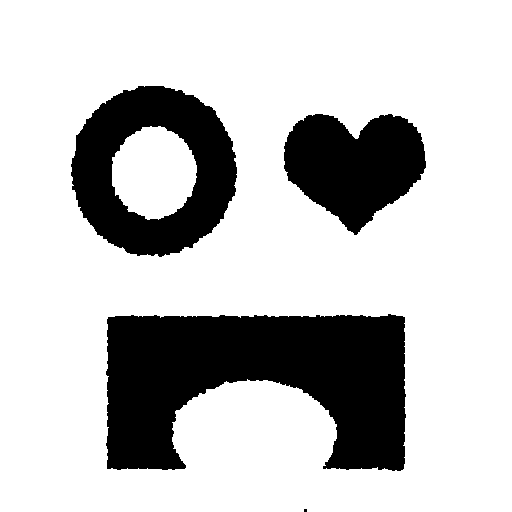}
        \caption{\centering $\sigma^2 = 0.5$}
    \end{subfigure}
    \caption{Results for noisy images using modified approach. \textit{Left column}: raw images, \textit{middle column}: filtered images, \textit{right column}: segmented images.}
    \label{fig5}
\end{figure}

\begin{figure}[h]
    \begin{subfigure}{.32 \textwidth}
        \centering
        \includegraphics[width=2.5cm]{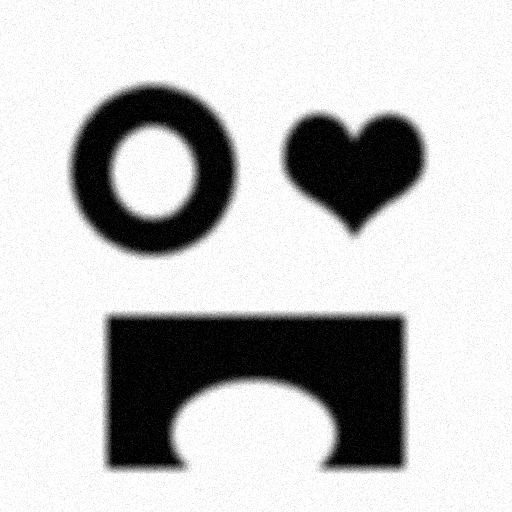}
		\caption{\centering $\mathcal{B}_g$}
    \end{subfigure}
    \begin{subfigure}{.32 \textwidth}
        \centering
        \includegraphics[width=2.5cm]{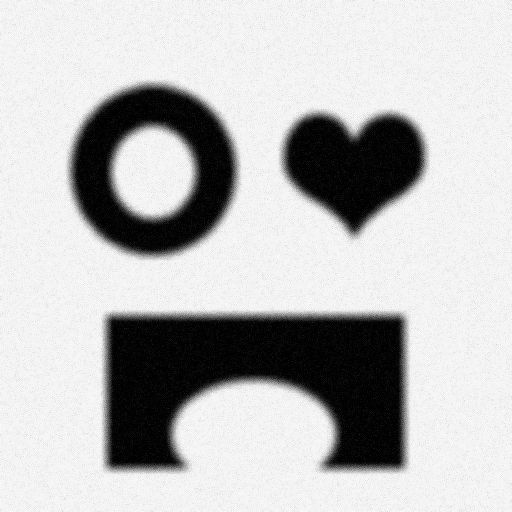}
		\caption{\centering $\mathcal{B}_g$}
    \end{subfigure}
    \begin{subfigure}{.32 \textwidth}
        \centering
        \includegraphics[width=2.5cm]{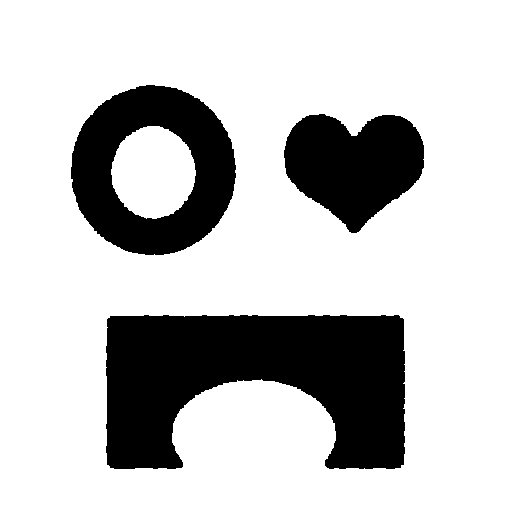}
		\caption{\centering $\mathcal{B}_g$}
    \end{subfigure}\\[0.1cm]
     \begin{subfigure}{.32 \textwidth}
        \centering
        \includegraphics[width=2.5cm]{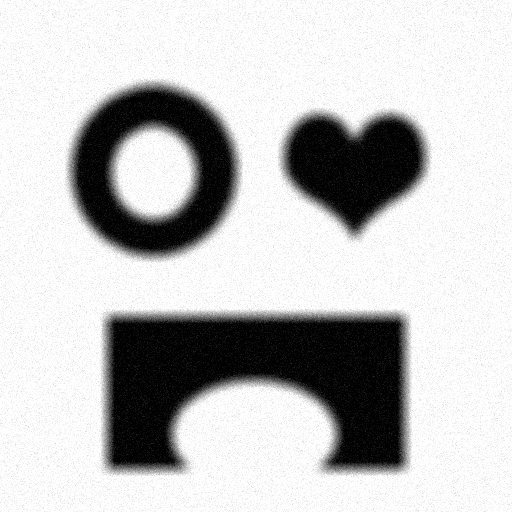}
		\caption{\centering $\mathcal{B}_a$}
    \end{subfigure}
    \begin{subfigure}{.32 \textwidth}
        \centering
        \includegraphics[width=2.5cm]{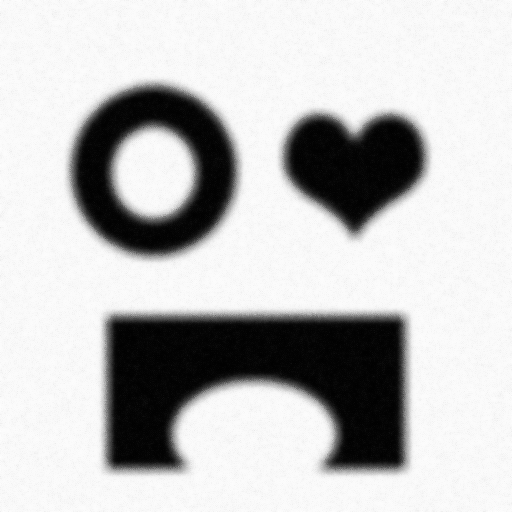}
		\caption{\centering $\mathcal{B}_a$}
    \end{subfigure}
    \begin{subfigure}{.32 \textwidth}
        \centering
        \includegraphics[width=2.5cm]{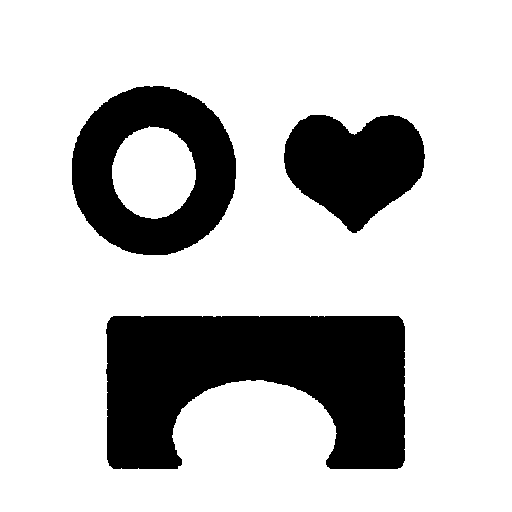}
		\caption{\centering $\mathcal{B}_a$}
    \end{subfigure}\\[0.1cm]
    \begin{subfigure}{.32 \textwidth}
        \centering
        \includegraphics[width=2.5cm]{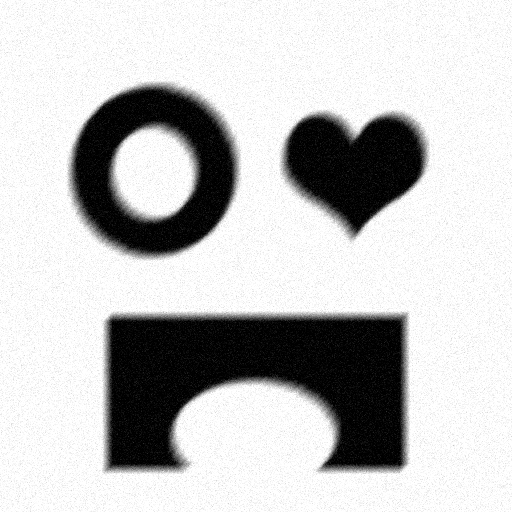}
		\caption{\centering $\mathcal{B}_m$}
    \end{subfigure}
    \begin{subfigure}{.32 \textwidth}
        \centering
        \includegraphics[width=2.5cm]{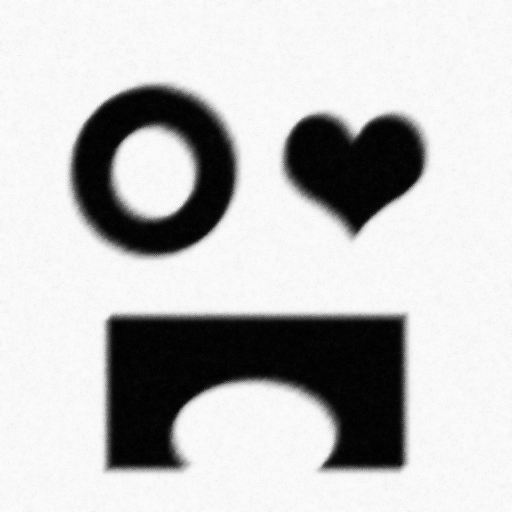}
		\caption{\centering $\mathcal{B}_m$}
    \end{subfigure}
    \begin{subfigure}{.32 \textwidth}
        \centering
        \includegraphics[width=2.5cm]{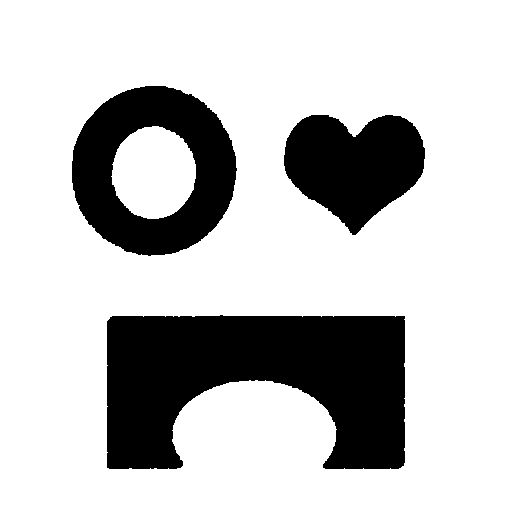}
		\caption{\centering $\mathcal{B}_m$}
    \end{subfigure}
    \caption{Results for blurry images using modified approach. \textit{Left column}: raw images, \textit{middle column}: filtered images, \textit{right column}: segmented images.}
    \label{fig6}
\end{figure}

In this subsection, we evaluate the performance of the modified approach. Filtering and segmentation results are presented in Fig. \ref{fig5} and Fig. \ref{fig6}, with a quantitative comparison provided in Tables \ref{noise2-compare} and \ref{blur2-compare}. Similar to the basic approach, the modified one demonstrates robust segmentation results. The quantitative analysis indicates that it achieves higher values of JS, DSC, and SA for the images with $\sigma^2 = 0.3$ and $\sigma^2 = 0.5$. A direct comparison between the basic and modified approaches, presented in Tables \ref{noise3-compare} and \ref{blur3-compare}, shows that the modified method generally yields slightly better metric values. However, for $\sigma^2 = 0.5$ and $\mb_g,~\mb_a$ the basic approach shows JS, DSC and SA scores that are up to $0.31 \%$ higher.

\begin{figure}[H]
    \begin{subfigure}{.32 \textwidth}
        \centering
        \includegraphics[width=2.5cm]{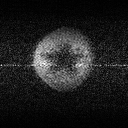}
    \end{subfigure}
    \begin{subfigure}{.32 \textwidth}
        \centering
        \includegraphics[width=2.5cm]{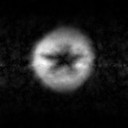}
    \end{subfigure}
    \begin{subfigure}{.32 \textwidth}
        \centering
        \includegraphics[width=2.5cm]{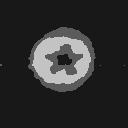}
    \end{subfigure}
    \caption{Results for the papaya MRI scan using modified approach. \textit{Left to right}: raw image, filtered image, segmented image.}
    \label{fig7}
\end{figure}

The results for the MRI scan of a papaya, shown in Fig. \ref{fig7}, further illustrate the robustness of the modified approach in filtering and segmenting low-field MRI scans.

All parameter values used in the experiments are provided in \ref{app1}. Additional tests with the total variation diffusion coefficient (Eq.~\ref{eq7}) were also conducted but did not yield improved results.

\renewcommand{\arraystretch}{1.5}
\begin{table}[H]
\captionsetup{justification=centering}
\caption{Modified approach, images with different levels of Gaussian noise added (in $\%$).}
\centering
\begin{tabular}{ m{2cm}  m{1.2cm} m{1.2cm}  m{1.2cm} m{1.2cm}  m{1.2cm} m{1.2cm}}
\hline
& \multicolumn{2}{c}{JS} & \multicolumn{2}{c}{DSC}&  \multicolumn{2}{c}{SA} \\
\hline
Image & TSIS & Ours & TSIS & Ours & TSIS & Ours\\
\hline
$\sigma^2 = 0.1$ & \textbf{98.9} & 98.55 & \textbf{99.4} & 99.27 & \textbf{99.7} & 99.66\\
$\sigma^2 = 0.3$ & 97.6 & \textbf{97.96} & 98.8 & \textbf{98.97} & 99.4 & \textbf{99.53}\\
$\sigma^2 = 0.5$ & 96.53 & \textbf{97.51} & 98.2 & \textbf{98.74} & 99.2 & \textbf{99.42}\\
\hline
\end{tabular}
\label{noise2-compare}
\end{table}

\renewcommand{\arraystretch}{1.5}
\begin{table}[H]
\captionsetup{justification=centering}
\caption{Modified approach, images with different types of blurring (in $\%$).}
\centering
\begin{tabular}{ m{2cm}  m{1.2cm} m{1.2cm}  m{1.2cm} m{1.2cm}  m{1.2cm} m{1.2cm}}
\hline
& \multicolumn{2}{c}{JS} & \multicolumn{2}{c}{DSC}&  \multicolumn{2}{c}{SA} \\
\hline
Image & TSIS & Ours & TSIS & Ours & TSIS & Ours\\
\hline
$\mathcal{B}_g$ & \textbf{98.9} & 97.38 & \textbf{98.9} & 98.67 & \textbf{99.5} & 99.39\\
$\mathcal{B}_a$ & \textbf{98.3} & 97.62 & 97.7 & \textbf{98.80} & 98.9 & \textbf{99.44}\\
$\mathcal{B}_m$ & \textbf{99.4} & 98.11 & \textbf{99.5} & 99.04 & \textbf{99.9} & 99.56\\
\hline
\end{tabular}
\label{blur2-compare}
\end{table}

\renewcommand{\arraystretch}{1.5}
\begin{table}[H]
\captionsetup{justification=centering}
\caption{Basic/modified approach, images with different levels of Gaussian noise added (in $\%$).}
\centering
\begin{tabular}{m{2cm} m{1.2cm} m{1.3cm}  m{1.2cm} m{1.3cm}  m{1.2cm} m{1.3cm}}
\hline
& \multicolumn{2}{c}{JS} & \multicolumn{2}{c}{DSC}&  \multicolumn{2}{c}{SA} \\
\hline
Image & Basic & Modified & Basic & Modified & Basic & Modified\\
\hline
$\sigma^2 = 0.1$ & 98.52 & \textbf{98.55} & 99.25 & \textbf{99.27} & \textbf{99.66} & \textbf{99.66}\\
$\sigma^2 = 0.3$ & 97.80 & \textbf{97.96} & 98.89 & \textbf{98.97} & 99.49 & \textbf{99.53}\\
$\sigma^2 = 0.5$ & \textbf{97.51} & 97.45 & \textbf{98.74} & 98.71 & \textbf{99.42} & 99.40\\
\hline
\end{tabular}
\label{noise3-compare}
\end{table}

\renewcommand{\arraystretch}{1.5}
\begin{table}[H]
\captionsetup{justification=centering}
\caption{Basic/modified approach, images with different types of blurring (in $\%$).}
\centering
\begin{tabular}{ m{2cm}  m{1.2cm} m{1.3cm}  m{1.2cm} m{1.3cm}  m{1.2cm} m{1.3cm}}
\hline
& \multicolumn{2}{c}{JS} & \multicolumn{2}{c}{DSC}&  \multicolumn{2}{c}{SA} \\
\hline
Image & Basic & Modified & Basic & Modified & Basic & Modified\\
\hline
$\mathcal{B}_g$ & \textbf{97.69} & 97.38 & 98.43 & \textbf{98.67} & \textbf{99.46} & 99.39\\
$\mathcal{B}_a$ & \textbf{97.74} & 97.62 & \textbf{98.86} & 98.80 & \textbf{99.47} & 99.44\\
$\mathcal{B}_m$ & 97.48 & \textbf{98.11} & 98.73 & \textbf{99.04} & 99.41 & \textbf{99.56}\\
\hline
\end{tabular}
\label{blur3-compare}
\end{table}

\section{Conclusions}
\label{S.5}

We have investigated a multi-stage PDE-based image processing approach consisting of three key stages: image reconstruction, denoising, and segmentation. In the first stage, the inverse Fourier transform was applied to convert MRI frequency-domain data into an initial spatial image. In the second stage, denoising techniques were applied to the image: for that, we employed the Nordstrom model \cite{nordstrom1990biased} alongside a generalized elastic net diffusion coefficient. The resulting system was discretised and solved using several numerical techniques described in \cite{shan2022deflated}. In the last stage, the final segmentation was obtained.  Additionally, we compared this approach to a modified version, in which the diffusion coefficient was computed from a presegmented image combined with Gaussian filtering. This modification aimed to produce sharper jumps in the values of the diffusion coefficient, particularly at object edges, enhancing segmentation accuracy.

From our analysis, we draw the following conclusions. Both the basic and modified image processing approaches demonstrated robust performance when ground truth (GT) segmentation was available, allowing for quantitative evaluation. This analysis showed that our approaches yielded higher evaluation metrics on noisier images compared to those reported in \cite{wu2021two}. However, the modified approach did not consistently outperform the basic one. This can be explained by the fact that the basic performance was already nearly perfect, with evaluation criteria reaching a minimum of approximately $97.5 \%$, leaving limited room for further improvement.

\section*{Acknowledgements}
Ksenia Slepova was supported by the Peter Paul Peterich Foundation via the TU Delft University Fund. The authors express their gratitude to the MRI Laboratory team at Mbarara University of Science and Technology, led by Dr. Johnes Obungoloch, for their hospitality and cooperation during the research visit to Uganda.

\appendix
\section{Parameter Settings for Experiments}
\label{app1}

In this appendix, we present the values of the parameters used in the experiments (see Table  \ref{parameters}). $K, d_p, \eta$ are the parameters involved in the diffusion coefficient computation (see section \ref{noise-reduction}), $n_{\text{cl.}}$ is the maximum number of clusters used in the segmentation, and $n_{\text{b}}$ is the pixel size of the patch employed in the operation of the morphological closing \cite{shapiro2001computer}.

\setcounter{table}{0}
\renewcommand{\arraystretch}{1.3}
\begin{table}[H]
\captionsetup{justification=centering}
\caption{Parameter settings for the experiments.}
\centering
\begin{tabular}{c | c | c | c | c | c | c}
Test image & Approach & $K$ & $d_p$ & $\eta$ & $n_{\text{cl.}}$ & $n_b$ \\
\hline
Synthetic $\sigma^2 = 0.1$ & Basic & $0.1$ & $1.2$ & $10^4$ & $2$ & $1$ \\
Synthetic $\sigma^2 = 0.1$ & Modified & $0.1$ & $1.9$ & $10^5$ & $2$ & $2$ \\
Synthetic $\sigma^2 = 0.3$ & Basic & $0.1$ & $1.8$ & $10^4$ & $2$ & $1$\\
Synthetic $\sigma^2 = 0.3$ & Modified &  $0.1$ & $1.9$ & $10^5$ & $2$ & $1$ \\
Synthetic $\sigma^2 = 0.5$ & Basic & $0.1$ & $1.8$ & $10^4$ & $2$ & $1$ \\
Synthetic $\sigma^2 = 0.5$ & Modified & $0.1$ & $1.9$ & $10^5$ & $2$ & $1$ \\
Synthetic $\mathcal{B}_g$ & Basic & $0.1$ & $1.0$ & $10^4$ & $2$ & $0$ \\
Synthetic $\mathcal{B}_g$ & Modified & $0.1$ & $1.8$ & $10^6$ & $2$ & $2$ \\
Synthetic $\mathcal{B}_a$ & Basic & $0.1$ & $1.0$ & $10^4$ & $2$ & $0$ \\
Synthetic $\mathcal{B}_a$ & Modified & $0.1$ & $1.5$ & $10^6$ & $2$ & $2$ \\
Synthetic $\mathcal{B}_m$ & Basic & $0.1$ & $1.0$ & $10^5$ & $2$ & $0$ \\
Synthetic $\mathcal{B}_m$ & Modified & $0.1$ & $1.5$ & $10^6$ & $2$ & $1$ \\
Papaya MRI scan & Basic & $0.1$ & $1.9$ & $10^3$ & $3$ & $1$ \\
Papaya MRI scan & Modified & $0.1$ & $1.0$ & $10^5$ & $3$ & $1$ 
\end{tabular}
\label{parameters}
\end{table}

Standard deviation of the Gaussian smoothing was taken as $\sigma_f = 1.0$ for all the tests.

\section{Gaussian Smoothing}
\label{app2}


The Gaussian smoothing is performed by convolving the Gaussian smoothing operator with the desired image. The goal of such a convolution is to blur the image and remove detail. This operator has the following form:

\begin{align*}
\mathcal{G}(x,y) &= \frac{1}{(\sqrt{2\pi} \sigma_f)^{\frac{n}{2}}} \cdot \exp \left( - \frac{x^2 + y^2}{2\sigma_f^2} \right),
\end{align*}
where $n$ is the dimension of the problem ($2D$) and $\sigma_f$ is the standard deviation that defines the level of smoothing.
%



\bibliographystyle{elsarticle-num-names} 
\bibliography{report.bib}

\begin{thebibliography}{33}
\expandafter\ifx\csname natexlab\endcsname\relax\def\natexlab#1{#1}\fi
\providecommand{\url}[1]{\texttt{#1}}
\providecommand{\href}[2]{#2}
\providecommand{\path}[1]{#1}
\providecommand{\DOIprefix}{doi:}
\providecommand{\ArXivprefix}{arXiv:}
\providecommand{\URLprefix}{URL: }
\providecommand{\Pubmedprefix}{pmid:}
\providecommand{\doi}[1]{\href{http://dx.doi.org/#1}{\path{#1}}}
\providecommand{\Pubmed}[1]{\href{pmid:#1}{\path{#1}}}
\providecommand{\bibinfo}[2]{#2}
\ifx\xfnm\relax \def\xfnm[#1]{\unskip,\space#1}\fi
\bibitem[{O'Reilly et~al.(2021)O'Reilly, Teeuwisse, de~Gans, Koolstra, and
  Webb}]{o2021vivo}
\bibinfo{author}{T.~O'Reilly}, \bibinfo{author}{W.~M. Teeuwisse},
  \bibinfo{author}{D.~de~Gans}, \bibinfo{author}{K.~Koolstra},
  \bibinfo{author}{A.~G. Webb},
\newblock \bibinfo{title}{In vivo {3D} brain and extremity {MRI} at 50 {mT}
  using a permanent magnet {Halbach} array},
\newblock \bibinfo{journal}{Magnetic resonance in medicine}
  \bibinfo{volume}{85} (\bibinfo{year}{2021}) \bibinfo{pages}{495--505}.
  \DOIprefix\doi{10.1002/mrm.28396}.
\bibitem[{O'Reilly et~al.(2019)O'Reilly, Teeuwisse, and Webb}]{o2019three}
\bibinfo{author}{T.~O'Reilly}, \bibinfo{author}{W.~Teeuwisse},
  \bibinfo{author}{A.~Webb},
\newblock \bibinfo{title}{Three-dimensional {MRI} in a homogenous 27 cm
  diameter bore {Halbach} array magnet},
\newblock \bibinfo{journal}{Journal of Magnetic Resonance}
  \bibinfo{volume}{307} (\bibinfo{year}{2019}) \bibinfo{pages}{106578}.
  \DOIprefix\doi{10.1016/j.jmr.2019.106578}.
\bibitem[{Obungoloch et~al.(2023)Obungoloch, Muhumuza, Teeuwisse, Harper,
  Etoku, Asiimwe, Tusiime, Gombya, Mugume, Namutebi
  et~al.}]{obungoloch2023site}
\bibinfo{author}{J.~Obungoloch}, \bibinfo{author}{I.~Muhumuza},
  \bibinfo{author}{W.~Teeuwisse}, \bibinfo{author}{J.~Harper},
  \bibinfo{author}{I.~Etoku}, \bibinfo{author}{R.~Asiimwe},
  \bibinfo{author}{P.~Tusiime}, \bibinfo{author}{G.~Gombya},
  \bibinfo{author}{C.~Mugume}, \bibinfo{author}{M.~H. Namutebi}, et~al.,
\newblock \bibinfo{title}{On-site construction of a point-of-care low-field
  {MRI} system in {Africa}},
\newblock \bibinfo{journal}{NMR in Biomedicine}  (\bibinfo{year}{2023})
  \bibinfo{pages}{e4917}. \DOIprefix\doi{10.1002/nbm.4917}.
\bibitem[{Shan and van Gijzen(2022)}]{shan2022deflated}
\bibinfo{author}{X.~Shan}, \bibinfo{author}{M.~B. van Gijzen},
\newblock \bibinfo{title}{Deflated preconditioned {Conjugate} {Gradient}
  methods for noise filtering of low-field {MR} images},
\newblock \bibinfo{journal}{Journal of Computational and Applied Mathematics}
  \bibinfo{volume}{400} (\bibinfo{year}{2022}) \bibinfo{pages}{113730}.
  \DOIprefix\doi{10.1016/j.cam.2021.113730}.
\bibitem[{Weickert(2001)}]{weickert2001efficient}
\bibinfo{author}{J.~Weickert},
\newblock \bibinfo{title}{Efficient image segmentation using partial
  differential equations and morphology},
\newblock \bibinfo{journal}{Pattern Recognition} \bibinfo{volume}{34}
  (\bibinfo{year}{2001}) \bibinfo{pages}{1813--1824}.
  \DOIprefix\doi{10.1016/S0031-3203(00)00109-6}.
\bibitem[{Chen et~al.(2000)Chen, Vemuri, and Wang}]{chen2000image}
\bibinfo{author}{Y.~Chen}, \bibinfo{author}{B.~C. Vemuri},
  \bibinfo{author}{L.~Wang},
\newblock \bibinfo{title}{Image denoising and segmentation via nonlinear
  diffusion},
\newblock \bibinfo{journal}{Computers \& Mathematics with Applications}
  \bibinfo{volume}{39} (\bibinfo{year}{2000}) \bibinfo{pages}{131--149}.
  \DOIprefix\doi{10.1016/S0898-1221(00)00050-X}.
\bibitem[{Karasev et~al.(2013)Karasev, Kolesov, Fritscher, Vela, Mitchell, and
  Tannenbaum}]{karasev2013interactive}
\bibinfo{author}{P.~Karasev}, \bibinfo{author}{I.~Kolesov},
  \bibinfo{author}{K.~Fritscher}, \bibinfo{author}{P.~Vela},
  \bibinfo{author}{P.~Mitchell}, \bibinfo{author}{A.~Tannenbaum},
\newblock \bibinfo{title}{Interactive medical image segmentation using {PDE}
  control of active contours},
\newblock \bibinfo{journal}{IEEE transactions on medical imaging}
  \bibinfo{volume}{32} (\bibinfo{year}{2013}) \bibinfo{pages}{2127--2139}.
  \DOIprefix\doi{10.1109/tmi.2013.2274734}.
\bibitem[{Kollem et~al.(2019)Kollem, Reddy, and Rao}]{kollem2019review}
\bibinfo{author}{S.~Kollem}, \bibinfo{author}{K.~R. Reddy},
  \bibinfo{author}{D.~S. Rao},
\newblock \bibinfo{title}{A review of image denoising and segmentation methods
  based on medical images},
\newblock \bibinfo{journal}{International Journal of Machine Learning and
  Computing} \bibinfo{volume}{9} (\bibinfo{year}{2019})
  \bibinfo{pages}{288--295}. \DOIprefix\doi{10.18178/ijmlc.2019.9.3.800}.
\bibitem[{Kollem et~al.(2021)Kollem, Reddy, and Rao}]{kollem2021}
\bibinfo{author}{S.~Kollem}, \bibinfo{author}{K.~R. Reddy},
  \bibinfo{author}{D.~S. Rao},
\newblock \bibinfo{title}{Improved partial differential equation-based total
  variation approach to non-subsampled contourlet transform for medical image
  denoising},
\newblock \bibinfo{journal}{Multimedia Tools and Applications}
  \bibinfo{volume}{80} (\bibinfo{year}{2021}) \bibinfo{pages}{2663--2689}.
  \DOIprefix\doi{10.1007/s11042-020-09745-1}.
\bibitem[{Riya et~al.(2021)Riya, Gupta, and Lamba}]{RIYA2021106}
\bibinfo{author}{Riya}, \bibinfo{author}{B.~Gupta}, \bibinfo{author}{S.~S.
  Lamba},
\newblock \bibinfo{title}{An efficient anisotropic diffusion model for image
  denoising with edge preservation},
\newblock \bibinfo{journal}{Computers \& Mathematics with Applications}
  \bibinfo{volume}{93} (\bibinfo{year}{2021}) \bibinfo{pages}{106--119}.
  \DOIprefix\doi{10.1016/j.camwa.2021.03.029}.
\bibitem[{Ranjbarzadeh et~al.(2021)Ranjbarzadeh, Bagherian~Kasgari,
  Jafarzadeh~Ghoushchi, Anari, Naseri, and Bendechache}]{ranjbarzadeh2021brain}
\bibinfo{author}{R.~Ranjbarzadeh}, \bibinfo{author}{A.~Bagherian~Kasgari},
  \bibinfo{author}{S.~Jafarzadeh~Ghoushchi}, \bibinfo{author}{S.~Anari},
  \bibinfo{author}{M.~Naseri}, \bibinfo{author}{M.~Bendechache},
\newblock \bibinfo{title}{Brain tumor segmentation based on deep learning and
  an attention mechanism using {MRI} multi-modalities brain images},
\newblock \bibinfo{journal}{Scientific Reports} \bibinfo{volume}{11}
  (\bibinfo{year}{2021}) \bibinfo{pages}{1--17}.
  \DOIprefix\doi{10.1038/s41598-021-90428-8}.
\bibitem[{Liu et~al.(2023)Liu, Tong, Chen, Jiang, Zhou, Zhang, Zhang, Jin, and
  Zhou}]{liu2023deep}
\bibinfo{author}{Z.~Liu}, \bibinfo{author}{L.~Tong}, \bibinfo{author}{L.~Chen},
  \bibinfo{author}{Z.~Jiang}, \bibinfo{author}{F.~Zhou},
  \bibinfo{author}{Q.~Zhang}, \bibinfo{author}{X.~Zhang},
  \bibinfo{author}{Y.~Jin}, \bibinfo{author}{H.~Zhou},
\newblock \bibinfo{title}{Deep learning based brain tumor segmentation: a
  survey},
\newblock \bibinfo{journal}{Complex \& intelligent systems} \bibinfo{volume}{9}
  (\bibinfo{year}{2023}) \bibinfo{pages}{1001--1026}.
  \DOIprefix\doi{10.1007/s40747-022-00815-5}.
\bibitem[{Guo et~al.(2019)Guo, Li, Huang, Guo, and Li}]{guo2019deep}
\bibinfo{author}{Z.~Guo}, \bibinfo{author}{X.~Li}, \bibinfo{author}{H.~Huang},
  \bibinfo{author}{N.~Guo}, \bibinfo{author}{Q.~Li},
\newblock \bibinfo{title}{Deep learning-based image segmentation on multimodal
  medical imaging},
\newblock \bibinfo{journal}{IEEE Transactions on Radiation and Plasma Medical
  Sciences} \bibinfo{volume}{3} (\bibinfo{year}{2019})
  \bibinfo{pages}{162--169}. \DOIprefix\doi{10.1109/TRPMS.2018.2890359}.
\bibitem[{Wu et~al.(2021)Wu, Shao, Gu, Ng, and Zeng}]{wu2021two}
\bibinfo{author}{T.~Wu}, \bibinfo{author}{J.~Shao}, \bibinfo{author}{X.~Gu},
  \bibinfo{author}{M.~K. Ng}, \bibinfo{author}{T.~Zeng},
\newblock \bibinfo{title}{Two-stage image segmentation based on nonconvex
  $\ell^2 - \ell^p$ approximation and thresholding},
\newblock \bibinfo{journal}{Applied Mathematics and Computation}
  \bibinfo{volume}{403} (\bibinfo{year}{2021}) \bibinfo{pages}{126168}.
  \DOIprefix\doi{10.1016/j.amc.2021.126168}.
\bibitem[{MacQueen et~al.(1967)}]{macqueen1967some}
\bibinfo{author}{J.~MacQueen}, et~al.,
\newblock \bibinfo{title}{Some methods for classification and analysis of
  multivariate observations},
\newblock in: \bibinfo{booktitle}{Proceedings of the fifth Berkeley symposium
  on mathematical statistics and probability}, volume~\bibinfo{volume}{1},
  \bibinfo{organization}{Oakland, CA, USA}, \bibinfo{year}{1967}, pp.
  \bibinfo{pages}{281--297}.
\bibitem[{Adams and Bischof(1994)}]{adams1994seeded}
\bibinfo{author}{R.~Adams}, \bibinfo{author}{L.~Bischof},
\newblock \bibinfo{title}{Seeded region growing},
\newblock \bibinfo{journal}{IEEE Transactions on pattern analysis and machine
  intelligence} \bibinfo{volume}{16} (\bibinfo{year}{1994})
  \bibinfo{pages}{641--647}. \DOIprefix\doi{10.1109/34.295913}.
\bibitem[{Stockman and Shapiro(2001)}]{shapiro2001computer}
\bibinfo{author}{G.~Stockman}, \bibinfo{author}{L.~G. Shapiro},
  \bibinfo{title}{Computer {Vision}}, \bibinfo{edition}{1st} ed.,
  \bibinfo{publisher}{Prentice Hall PTR}, \bibinfo{address}{USA},
  \bibinfo{year}{2001}.
\bibitem[{Gonzalez and Woods(2006)}]{gonzalez2006dip}
\bibinfo{author}{R.~C. Gonzalez}, \bibinfo{author}{R.~E. Woods},
  \bibinfo{title}{Digital Image Processing (3rd Edition)},
  \bibinfo{publisher}{Prentice-Hall, Inc.}, \bibinfo{address}{USA},
  \bibinfo{year}{2006}.
\bibitem[{Jenks(1967)}]{jenks1967data}
\bibinfo{author}{G.~F. Jenks},
\newblock \bibinfo{title}{The data model concept in statistical mapping},
\newblock \bibinfo{journal}{International yearbook of cartography}
  \bibinfo{volume}{7} (\bibinfo{year}{1967}) \bibinfo{pages}{186--190}.
\bibitem[{Perona and Malik(1990)}]{perona1990scale}
\bibinfo{author}{P.~Perona}, \bibinfo{author}{J.~Malik},
\newblock \bibinfo{title}{Scale-space and edge detection using anisotropic
  diffusion},
\newblock \bibinfo{journal}{IEEE Transactions on pattern analysis and machine
  intelligence} \bibinfo{volume}{12} (\bibinfo{year}{1990})
  \bibinfo{pages}{629--639}. \DOIprefix\doi{10.1109/34.56205}.
\bibitem[{Weickert et~al.(1998)}]{weickert1998anisotropic}
\bibinfo{author}{J.~Weickert}, et~al., \bibinfo{title}{Anisotropic diffusion in
  image processing}, volume~\bibinfo{volume}{1}, \bibinfo{publisher}{Teubner
  Stuttgart}, \bibinfo{year}{1998}.
\bibitem[{Nordstr{\"o}m(1990)}]{nordstrom1990biased}
\bibinfo{author}{K.~N. Nordstr{\"o}m},
\newblock \bibinfo{title}{Biased anisotropic diffusion: a unified
  regularization and diffusion approach to edge detection},
\newblock \bibinfo{journal}{Image and vision computing} \bibinfo{volume}{8}
  (\bibinfo{year}{1990}) \bibinfo{pages}{318--327}.
  \DOIprefix\doi{10.1016/0262-8856(90)80008-H}.
\bibitem[{Whitaker and Pizer(1993)}]{whitaker1993multi}
\bibinfo{author}{R.~T. Whitaker}, \bibinfo{author}{S.~M. Pizer},
\newblock \bibinfo{title}{A multi-scale approach to nonuniform diffusion},
\newblock \bibinfo{journal}{CVGIP: image understanding} \bibinfo{volume}{57}
  (\bibinfo{year}{1993}) \bibinfo{pages}{99--110}.
  \DOIprefix\doi{10.1006/ciun.1993.1006}.
\bibitem[{You et~al.(1994)You, Kaveh, Xu, and Tannenbaum}]{you1994analysis}
\bibinfo{author}{Y.-L. You}, \bibinfo{author}{M.~Kaveh}, \bibinfo{author}{W.-Y.
  Xu}, \bibinfo{author}{A.~Tannenbaum},
\newblock \bibinfo{title}{Analysis and design of anisotropic diffusion for
  image processing},
\newblock in: \bibinfo{booktitle}{Proceedings of 1st International Conference
  on Image Processing}, volume~\bibinfo{volume}{2},
  \bibinfo{organization}{IEEE}, \bibinfo{year}{1994}, pp.
  \bibinfo{pages}{497--501}. \DOIprefix\doi{10.1109/ICIP.1994.413620}.
\bibitem[{Zou and Hastie(2005)}]{elasticnet}
\bibinfo{author}{H.~Zou}, \bibinfo{author}{T.~Hastie},
\newblock \bibinfo{title}{Regularization and {Variable} {Selection} {Via} the
  {Elastic} {Net}},
\newblock \bibinfo{journal}{Journal of the Royal Statistical Society Series B:
  Statistical Methodology} \bibinfo{volume}{67} (\bibinfo{year}{2005})
  \bibinfo{pages}{301--320}. \DOIprefix\doi{10.1111/j.1467-9868.2005.00503.x}.
\bibitem[{Vese and Le~Guyader(2015)}]{vese2015variational}
\bibinfo{author}{L.~A. Vese}, \bibinfo{author}{C.~Le~Guyader},
  \bibinfo{title}{Variational methods in image processing},
  \bibinfo{publisher}{CRC Press}, \bibinfo{year}{2015}.
  \DOIprefix\doi{10.1007/978-1-4684-0567-5}.
\bibitem[{Catt{\'e} et~al.(1992)Catt{\'e}, Lions, Morel, and
  Coll}]{catte1992image}
\bibinfo{author}{F.~Catt{\'e}}, \bibinfo{author}{P.-L. Lions},
  \bibinfo{author}{J.-M. Morel}, \bibinfo{author}{T.~Coll},
\newblock \bibinfo{title}{Image selective smoothing and edge detection by
  nonlinear diffusion},
\newblock \bibinfo{journal}{SIAM Journal on Numerical analysis}
  \bibinfo{volume}{29} (\bibinfo{year}{1992}) \bibinfo{pages}{182--193}.
  \DOIprefix\doi{10.1137/0729012}.
\bibitem[{Vogel and Oman(1996)}]{vogel1996total}
\bibinfo{author}{C.~R. Vogel}, \bibinfo{author}{M.~E. Oman},
\newblock \bibinfo{title}{Iterative {Methods} for {Total} {Variation}
  {Denoising}},
\newblock \bibinfo{journal}{SIAM Journal on Scientific Computing}
  \bibinfo{volume}{17} (\bibinfo{year}{1996}) \bibinfo{pages}{227--238}.
  \DOIprefix\doi{10.1137/0917016}.
\bibitem[{Nicolaides(1987)}]{nicolaides1987deflation}
\bibinfo{author}{R.~A. Nicolaides},
\newblock \bibinfo{title}{Deflation of conjugate gradients with applications to
  boundary value problems},
\newblock \bibinfo{journal}{SIAM Journal on Numerical Analysis}
  \bibinfo{volume}{24} (\bibinfo{year}{1987}) \bibinfo{pages}{355--365}.
  \DOIprefix\doi{10.1137/0724027}.
\bibitem[{Gallagher et~al.(2008)Gallagher, Nemeth, and
  Hacein-Bey}]{gallagher2008introduction}
\bibinfo{author}{T.~A. Gallagher}, \bibinfo{author}{A.~J. Nemeth},
  \bibinfo{author}{L.~Hacein-Bey},
\newblock \bibinfo{title}{An introduction to the {Fourier} transform:
  relationship to {MRI}},
\newblock \bibinfo{journal}{American Journal of Roentgenology}
  \bibinfo{volume}{190} (\bibinfo{year}{2008}) \bibinfo{pages}{1396--1405}.
  \DOIprefix\doi{10.2214/ajr.07.2874}.
\bibitem[{de~Leeuw~den Bouter(2022)}]{de2022image}
\bibinfo{author}{M.~de~Leeuw~den Bouter}, \bibinfo{title}{Image
  {Reconstruction} for {Low-Field} {MRI}}, Ph.D. thesis, Delft University of
  Technology, \bibinfo{year}{2022}.
  \DOIprefix\doi{10.4233/uuid:f3c4431d-368c-4a17-aacf-2d1283688a1a}.
\bibitem[{Boyd et~al.(2011)Boyd, Parikh, Chu, Peleato, Eckstein
  et~al.}]{boyd2011distributed}
\bibinfo{author}{S.~Boyd}, \bibinfo{author}{N.~Parikh},
  \bibinfo{author}{E.~Chu}, \bibinfo{author}{B.~Peleato},
  \bibinfo{author}{J.~Eckstein}, et~al.,
\newblock \bibinfo{title}{Distributed optimization and statistical learning via
  the alternating direction method of multipliers},
\newblock \bibinfo{journal}{Foundations and Trends in Machine learning}
  \bibinfo{volume}{3} (\bibinfo{year}{2011}) \bibinfo{pages}{1--122}.
  \DOIprefix\doi{10.1561/2200000016}.
\bibitem[{M{\"u}ller et~al.(2022)M{\"u}ller, Soto-Rey, and
  Kramer}]{muller2022towards}
\bibinfo{author}{D.~M{\"u}ller}, \bibinfo{author}{I.~Soto-Rey},
  \bibinfo{author}{F.~Kramer},
\newblock \bibinfo{title}{Towards a guideline for evaluation metrics in medical
  image segmentation},
\newblock \bibinfo{journal}{BMC Research Notes} \bibinfo{volume}{15}
  (\bibinfo{year}{2022}) \bibinfo{pages}{210}.
  \DOIprefix\doi{10.1186/s13104-022-06096-y}.

\end{thebibliography}



\end{document}